\documentclass[10pt,reqno]{amsart}
\usepackage{amsmath,amsopn,amssymb,amsthm,multicol, dsfont, bbm}
\usepackage{hyperref}
\usepackage{color}
\usepackage{graphicx}

\hypersetup{
  colorlinks   = true, 
  urlcolor     = red, 
  linkcolor    = red, 
  citecolor   = blue 
}

\DeclareMathOperator{\Ad}{Ad}
\DeclareMathOperator{\ad}{ad}
\DeclareMathOperator{\Aut}{Aut}

\DeclareMathOperator{\tr}{tr}

\DeclareMathOperator{\Ric}{Ric}

\DeclareMathOperator{\Span}{span}

\newcommand{\fr}{\mathfrak}
\newcommand{\al}{\alpha}

\newcommand{\bb}{\mathbb}

\DeclareMathOperator{\SO}{SO}
\DeclareMathOperator{\s}{S}

\DeclareMathOperator{\Sp}{Sp}
\DeclareMathOperator{\SU}{SU}
\DeclareMathOperator{\U}{U}
\DeclareMathOperator{\G}{G}
\DeclareMathOperator{\F}{F}
\DeclareMathOperator{\E}{E}
\newcommand{\thickline}{\noalign{\hrule height 1pt}}

 \newtheorem{lemma} {Lemma} [section]
\newtheorem{theorem}[lemma]{Theorem} 
\newtheorem{remark}[lemma] {Remark} 
\newtheorem{prop} [lemma]{Proposition}  
\newtheorem{definition}[lemma] {Definition} 
 
\newtheorem{example}[lemma] {Example}

\usepackage{lscape}

\begin{document}

\title{Equigeodesics on some classes of homogeneous spaces} 
\author{Marina Statha}
\address{University of Thessaly, Department of Mathematics, GR-35100 Lamia and University of Patras, Department of Mathematics, GR-26500 Rion,  Greece}
\email{marinastatha@uth.gr} 
\medskip

\begin{abstract}

We study homogeneous curves on some classes of reductive homogeneous spaces $G/H$ which are geodesics with respect to any $G$-invariant metric on $G/H$.  These curves are called equigeodesics.  The spaces we consider are certain Stiefel manifolds $V_k\bb{R}^n$, generalized Wallach spaces and spheres. 
We give a characterization for algebraic equigeodesics on $V_2\bb{R}^n$, $V_4\bb{R}^6$, $\SO(6)/\SO(3)\cdot\SO(2)$, ${W}^{6} = \U(3)/\U(1)^3$, ${W}^{12} = \Sp(3)/\Sp(1)^3$, $\bb{S}^{2n+1}\cong\U(n+1)/\U(n)$ and $\bb{S}^{4n+3}\cong\Sp(n+1)/\Sp(n)$.


\medskip
\noindent 2020 {\it Mathematics Subject Classification.} Primary 53C25; Secondary 53C30.

\medskip
\noindent {\it Keywords}:  Homogeneous space, equigeodesic, equigeodesic vector, Stiefel manifold, generalized Wallach space, sphere.
\end{abstract}

\maketitle
 

\section{Introduction}
\markboth{Marina Statha}{Equigeodesics on some classes of homogeneous spaces}

\subsection{Introduction}
Let $(G/H, g)$ be a Riemannian homogeneous space.  A geodesic $\gamma(t)$ through the origin $o=eH$ is called a {\it homogeneous geodesic} if it is an orbit of a one-parameter subgroup of $G$, i.e.
$
\gamma(t)=(\exp tX)\cdot o,
$
where $\exp : \fr{g} \to G$ is the exponential map and $X$ is a non zero vector in the Lie algebra $\fr{g}$ of $G$.  If all geodesics on $G/H$ are homogeneous geodesics the homogeneous space is called a {\it g.o. space} (from ``geodesic ordit").  The terminology was introduced by O. Kowalski and L. Vanhecke in \cite{KoVa}, who initiated a systematic study of such spaces.  Examples of such spaces are the symmetric spaces, isotropy irreducible homogeneous spaces, the weakly symmetric spaces and the naturally reductive spaces. A Riemannian homogeneous space $(G/H, g)$ is called {\it naturally reductive} if there exists a subspace $\fr{m}$ of $\fr{g}$ with $\fr{g} = \fr{h}\oplus\fr{m}$ and $\Ad^{H}\fr{m}\subset\fr{m}$, such that the endomorphism $\ad(X) :\fr{m}\to \fr{m}$ is skew-symmetric with respect to some $\Ad^{H}$-invariant inner product on $\fr{m}$, for any $X\in\fr{h}$.  At this point it is worth mentioning that the class of g.o spaces is larger than the class of naturally reductive spaces. In fact, A. Kaplan (\cite{Ka}) gave the first example of a g.o. space which is not naturally reductive. 

In \cite{CGN} the authors introduced the notion of {\it homogeneous equigeodesics}.  In fact, they studied homogeneous curves on generalized flag manifolds of type $A_{\ell}$, that are geodesics with respect to any $G$-invariant metric.  The notion of a equigeodesic is important to mechanics, since homogeneous geodesics describe relative equilibria of a dynamical system represented by an invariant metric and we consider the case when this property is stable under change of the metric.
 Later in \cite{WaZh}, Y. Wang and G. Zhao used the work \cite{GrNe} to describe equigeodesics on generalized flag manifolds with two isotropy summands.  They gave a general formula for finding equigeodesic vectors on generalized flag manifolds with second Betti number equal to one (that is flag manifolds which are determined by painting one black node in their Dynkin diagram).  Recently, the author studied equigeodesics on generalized flag manifolds with $\G_2$-type $\fr{t}$-roots (\cite{Sta2}).  The approach was based on the fact that the infinitesimal generator of the one parameter subgroup is an element of the Lie algebra of $G$.  So, it is natural to characterize the equigeodesics in terms of their infinitesimal generator.  This allows us to use a Lie theoretical approach for the study of homogeneous geodesics on flag manifolds.  The infinitesimal generator of an equigeodesic is called {\it equigeodesic vector}.  Actually, a vector is equigeodesic if and only if it is a solution of an algebraic nonlinear system of equations (see Proposition \ref{sxesi}), whose variables are the coefficients of the vector. However, there exist subspaces of the tangent space (at the origin) of homogeneous space $G/H$, all of whose elements are equigeodesic vectors.  The existence of these subspace depends on the geometric structure of $G/H$.  The equigeodesic vectors are separated into two classes: the {\it structural} and the {\it algebraic}.  For flag manifolds it is easier to study structural equigeodesic vectors (\cite{GrNe}, \cite{WaZh}, \cite{Sta2}), since the geometric structure can be expressed in terms of the Lie algebra structure of Lie group and more precisely the root space decomposition.    


\subsection{Outline}
In the present work we describe algebraic equigeodesic vectors on certain reductive homogeneous spaces such as Stiefel manifolds, generalized Wallach spaces and spheres.  Also, we introduce the {\it homogeneous Einstein equigeodesics}, that is homogeneous equigeodesics which are geodesic with respect to every Einstein metric of the homogeneous space.  Recall that an Einstein metric on a Riemannian manifold $(M, g)$ is the solution of the Einstein equation $\Ric_{g} = cg$, where $\Ric_g$ is the Ricci tensor corresponding to $g$ and $c$ a real number (called Einstein constant).  Besides the detailed exposition on Einstein manifolds in \cite{Be}, we refer to \cite{Wa1}, \cite{Wa2}, \cite{Arv3}, for more recent results. 
The main results of the paper are the following:

\smallskip
\noindent
{\bf Proposition A.}
(1) The generalized Wallach space $\SO(6)/\SO(3)\cdot\SO(2)$ admits six classes of algebraic equigeodesics vectors.  One of these is trivial.

\noindent
(2) The Wallach space $W^{6} = \U(3)/\U(1)^3$ admits only trivial equigeodesic vectors and the Wallach space $W^{12} = \Sp(3)/\Sp(1)^{3}$ 
admits five classes of algebraic equigeodesic vectors.  Four of these are trivial.

\smallskip
\noindent
{\bf Proposition B.}
The Stiefel manifold $V_{2}\bb{R}^n$ admits $n$ classes of algebraic equigeodesic vectors.  Two of these classes are trivial.  Also, the Einstein manifold $(V_{2}\bb{R}^n, \Lambda = (1,\lambda,\lambda))$ admits one class of trivial Einstein equigeodesic vectors and one class of algebraic Einstein equigeodesic vectors.

\smallskip
\noindent
{\bf Proposition C.}
The Stiefel manifold $V_{4}\bb{R}^{6} \to \SO(6)/\SO(3)\cdot\SO(2)$ admits, with respect to  Jensen's Einstein metrics, $42$ classes of algebraic equigeodesic vectors.

\smallskip
\noindent
{\bf Proposition D.}
The spheres $\bb{S}^{2n+1}\cong\U(n+1)/\U(n)$ and $\bb{S}^{4n+3}\cong\Sp(n+1)/\Sp(n)$ admit only trivial equigeodesic vectors.


\subsection{$G$-invariant metrics on reductive homogeneous spaces}\label{subsection1.3}
Let $G$ be a compact Lie group and $H$ be a closed subgroup so that $G$ acts almost effectively on $G/H$.  Let $\fr{g},$ $\fr{h}$ be the Lie algebras of $G$ and $H$, and $\Ad^{G} : G \to \Aut(\fr{g})$ be the adjoint representation of $G$.  A homogeneous space $G/H$ is called {\it reductive} if there is a decomposition $\fr{g} = \fr{h}\oplus\fr{m}$, which satisfies $\Ad(h)\fr{m}\subset\fr{m}$ for all $h\in H$.  In general, the reductive decomposition $\fr{m}$ may not be unique and a sufficient condition of existence is the case where the subgroup $H$ is compact.  It is well known that a Riemannian homogeneous space always admits a reductive decomposition (e.g. \cite{KoSz}).   For reductive homogeneous spaces there is always a natural identification of $\fr{m}$ with the tangent space $T_{eH}(G/H)$, given by $X\leftrightarrow X^{*}$ where $X^{*} = \frac{d}{dt}(\exp tX)\cdot eH\big|_{t = 0}$.  If $G$ is in addition semisimple, the negative of the Killing form $B$ of $\fr{g}$ is an $\Ad^{G}$-invariant inner product on $\fr{g}$, therefore we can choose the above decomposition with respect to this form.  A Riemannian metric $g$ on a homogeneous space $G/H$ is called $G$-invariant if the diffeomorphism $\tau_{\al} : G/H \to G/H,$ $\tau_{\al}(b H) = \al bH$ is an isometry, that is $\tau_{\al}^{*}g = g$.  In order to describe all $G$-invariant metrics on $G/H$ we need to compute the isotropy representation $\chi : H \to \Aut(\fr{m})$, $\chi(h) := (d\tau_h)_{eH} : \fr{m} \to \fr{m}$.  The next proposition is useful to compute the isotropy representation of a reductive homogeneous space. 

\begin{prop}\textnormal{(\cite{Arv1})}\label{isotrepr}
Let $G/H$ be a reductive homogeneous space and let $\fr{g} = \fr{h}\oplus\fr{m}$ be a reductive decomposition of $\fr{g}$.  Let $h\in H$, $X\in \fr{h}$ and $Y\in\fr{m}$.  Then
$
\Ad^{G}(h)(X + Y) = \Ad^{H}(h)X + \chi(h)Y
$
that is, the restriction $\Ad^{G}|_{H}$ splits into the sum $\Ad^{H}\oplus\chi$.  
\end{prop}  

Any $G$-invariant Riemannian metric on $G/H$ is in one-to-one correspondence with $\Ad^{H}$-invariant inner products $\left\langle \cdot, \cdot\right\rangle$ on $\fr{m}$, i.e. $\langle \Ad(h)X, \Ad(h)Y\rangle = \langle X, Y\rangle$ for all $h\in H$ and $X, Y\in \fr{m}$.  If $H$ is compact then $\fr{m} = \fr{h}^{\perp}$ with respect to $B$, the negative of the Killing form of $\fr{g}$.  Then the $G$-invariant metrics on $G/H$ are in one-to-one correspondence with operators $\Lambda : \fr{m}\to \fr{m}$, which are  symmetric with respect to $B$, positive definite and $\Ad^H$-equivariant, that is $\Ad(h)\circ\Lambda = \Lambda\circ\Ad(h)$ for all $h\in H$.
This relation is given by
$
\langle X, Y\rangle = B(\Lambda X, Y).
$
 Moreover, if $\fr{m}$ decomposes into a direct sum of $\Ad^{H}$-invariant irreducible and pairwise inequivalent modules $\fr{m}_{i}$, 
that is 
$
\fr{m} = \fr{m}_{1}\oplus\cdots\oplus\fr{m}_{s}, 
$
then all $\Ad^{H}$-invariant inner products on $\fr{m}$ are given by
$
\langle\cdot ,  \cdot\rangle = \lambda_{1}B|_{\fr{m}_{1}} + \cdots + \lambda_{s}B|_{\fr{m}_{s}}, \ \lambda_{i}\in\bb{R}^{+}, \ i = 1,\ldots,s.
$
For such inner products and $G$-invariant metrics we may use the common notation 
 $\Lambda = \lambda_1{\rm Id}_{\fr{m}_1} + \cdots + \lambda_{s}{\rm Id}_{\fr{m}_s}$.

\begin{definition}
Let $(M, g)$ a Riemannian manifold.  Then $M$ is called Einstein if the metric $g$ satisfies the equation $\Ric_g = c g$, where $\Ric_g$ is the Ricci tensor of $g$ and $c$ is a real number.
\end{definition}

\section{Generalized Wallach spaces, Stiefel manifold and spheres}

\subsection{Generalized Wallach spaces}\label{subsectionWallach}
Let $G/H$ be a reductive homogeneous space with $G$ a compact and semisimple Lie group and $H$ a compact subgroup of $G$.  Let $\fr{g} = \fr{h}\oplus\fr{m}$ be the reductive decomposition of $G/H$.  Then $G/H$ is called a {\it generalized Wallach space} if the module $\fr{m}$ decomposes into a direct sum of three $\Ad^{H}$-invariant irreducible modules pairwise orthogonal with respect to $B$, i.e. $\fr{m} = \fr{m}_{12}\oplus\fr{m}_{13}\oplus\fr{m}_{23}$, such that $[\fr{m}_{ij}, \fr{m}_{ij}]\subset\fr{h}$, for $1\leq i < j \leq 3$.  Every generalized Wallach space admits a three parameter family of $G$-invariant Riemannian metrics determined by $\Ad^{H}$-invariant inner products
\begin{equation}\label{Nikonorov}
\Lambda_{{\tiny\mathsf{Wal}}} = \lambda_{12}{\rm Id}_{\fr{m}_{12}} + \lambda_{13}{\rm Id}_{\fr{m}_{13}} + \lambda_{23}{\rm Id}_{\fr{m}_{23}},
\end{equation}  
where $\lambda_{12}, \lambda_{13},\lambda_{23}$ are positive real numbers.  

The classification of generalized Wallach spaces $G/H$ was obtained in \cite{Ni1} and \cite{ChKaLi}.  Some examples of such spaces are the following:

$\bullet$ {\it The Wallach spaces:}
$
{W}^6 = \SU(3)/{\rm T}_{\mathsf{max}},\ {W}^{12} = \Sp(3)/\Sp(1)^3,\ {W}^{24} = \F_4/{\rm Spin}(8)
$

$\bullet$ {\it The generalized flag manifolds:}

\centerline{$
\SU(\ell+m+n)/\s(\U(\ell)\cdot\U(m)\cdot\U(n)),\  \SO(2\ell)/\U(1)\cdot\U(\ell-1),\  \E_6/\U(1)\cdot\U(1)\cdot{\rm Spin}(8)
$}

$\bullet$ {\it The homogeneous spaces:}

\centerline{$
\SO(n_1+n_2+n_3)/\SO(n_1)\cdot\SO(n_2)\cdot\SO(n_3),\  \Sp(n_1+n_2+n_3)/\Sp(n_1)\cdot\Sp(n_2)\cdot\Sp(n_3)
$}

\noindent
The modules $\fr{m}_{12}, \fr{m}_{13}, \fr{m}_{23}$ satisfy the following property, stated in \cite{Ni2} and a proof can be found in \cite{ArSo}.
\begin{lemma}\label{BrackW}
Let $G/H$ a generalized Wallach space.  Then the $\Ad^{H}$-invariant summands $\fr{m}_{ij}$, $1\leq i<j \leq 3$ satisfy the Lie bracket relations
$
[\fr{m}_{12}, \fr{m}_{13}] \subset \fr{m}_{23}, \, [\fr{m}_{12}, \fr{m}_{23}]\subset \fr{m}_{13}, \, [\fr{m}_{13}, \fr{m}_{23}]\subset \fr{m}_{12}.
$
\end{lemma}

In Section \ref{section3} we study equigeodesics in the generalized Wallach spaces {\bf (I)} $\SO(n_1+n_2+n_3)/\SO(n_1)\cdot\SO(n_2)\cdot\SO(n_3)$, {\bf (II)} $W^6 = \U(3)/\U(1)^3$ and {\bf (III)} $W^{12} = \Sp(3)/\Sp(1)^3$, so for these spaces we explicitly describe below the decomposition $\fr{m} = \fr{m}_{12}\oplus\fr{m}_{13}\oplus\fr{m}_{23}$.

\smallskip
\noindent \underline{{\bf Case} {\bf{(I)}}} $G/H = \SO(n)/\SO(n_1)\cdot\SO(n_2)\cdot\SO(n_3), n=n_1+n_2+n_3.$
We embed the group $H = \SO(n_1)\cdot\SO(n_2)\cdot\SO(n_3)$ diagonally in $G = \SO(n_1 + n_2 + n_3)$ ($n_1\ge n_2\ge n_3\ge 2$). Then the tangent space $\fr{m}$ of $G/H$ is given by $\fr{h}^{\perp}$ in $ \fr{g} = \fr{so}(k_1+ k_2+k_3)$ with respect to  $B : \fr{so}(n)\times\fr{so}(n)\to\bb{R}$,  given by $B(X, Y) = -{\rm tr}XY$.  
Then we see that 
  $ \fr{m}$ is given by 
  \begin{equation*}
\fr{m}=  \left\{\begin{pmatrix}
 0 & {\fr{m}}_{12} & {\fr{m}}_{13}\\
 -{}^{t}_{}\!{\fr{m}}_{12} & 0 & {\fr{m}}_{23}\\
 -{}^{t}_{}\!{\fr{m}}_{13} & -{}^{t}_{}\!{\fr{m}}_{23} & 0 
 \end{pmatrix} \  \Big\vert \ 
 \begin{array}{l}
  \dim{\fr{m}}_{12} =n_1n_2,\\ \dim{\fr{m}}_{13}=n_1n_3,
  \\ \dim{\fr{m}}_{23}=n_2n_3 
  \end{array}
   \right\}. 
 \end{equation*}

Let $E_{ab}$ be the $n\times n$ matrix with $1$ at the $(ab)$-entry and $0$ elsewhere.  Then the set
$\mathcal{B}=\{\xi_{ab}=E_{ab}-E_{ba}: 1\le a<b\le n\}$ constitutes a $B$-orthogonal basis of $\fr{so}(n)$.
Note that $\xi_{ba}=-\xi_{ab}$, thus we have the following:

\begin{lemma}\label{bracSO}
If all four indices are distinct, then the Lie brackets in $\mathcal{B}$ are zero.
Otherwise,
$[\xi_{ab}, \xi_{bc}]=\xi_{ac}$, where $a,b,c$ are distinct.
\end{lemma}

\begin{remark}
\textnormal{By using the above Lemma it is easy to see that in the case where $n_1 = n_2 = n_3 = 2$, or $n_{i} = n_{j} = 2$ for $1\leq i < j\leq 3$, the homogeneous space $G/H=\SO(n_1+n_2+n_3)/\SO(n_1)\cdot\SO(n_2)\cdot\SO(n_3)$ it is not a generalized Wallach space.  In the first case the tangent space decomposes into six $\Ad^H$-invariant, irreducible and non equivalent summands and for the second case the tangent space decomposes into four $\Ad^{H}$-invariant, irreducible and non equivalent isotropy summands.  More precisely, if two of $n_i$'s are equal to 2, the two of these four summands of the tangent space of $G/H$ is identified to the tangent space of the Grassmannian ${\rm Gr}_{2}\bb{R}^4\cong\SO(4)/\SO(2)\cdot\SO(2) = \bb{S}^2\times\bb{S}^2$.  
}
\end{remark}

\noindent \underline{{\bf Case} {\bf (II)}} $G/H = \U(3)/\U(1)^3$.  We consider the embedding $\fr{h} = \fr{u}(1)\oplus\fr{u}(1)\oplus\fr{u}(1) \hookrightarrow 
\begin{pmatrix}
ia & 0 & 0\\
0& ib & 0\\
0 & 0 & ic
\end{pmatrix} \in\fr{u}(3)
$  
Then, the tangent space $\fr{m} \cong T_{eH}(G/H)$ is given by $\fr{h}^{\perp}$ in $\fr{u}(3)$ with respect to $B : \fr{u}(3)\times\fr{u}(3) \to \bb{R}$, $(X, Y) \mapsto -\tr{XY}$.  It is easy to see that $\fr{m}$ has the following form:
\begin{equation*}
 \left\{\begin{pmatrix}
 0 & {a}_{12} & {a}_{13}\\
 -{}^{}_{}\!{\bar{a}}_{12} & 0 & {a}_{23}\\
 -{}^{}_{}\!{\bar{a}}_{13} & -{}^{}_{}\!{\bar{a}}_{23} & 0 
 \end{pmatrix} \  \Big\vert \ 
  a_{ij}\in\bb{C},\ 1\leq i < j \leq 3
   \right\}.
\end{equation*}
  We set
\begin{eqnarray}\label{bracetSU}
e_{ij} = E_{ij}-E_{ji}, \ f_{ij} = \sqrt{-1}(E_{ij} + E_{ji}).
\end{eqnarray}
Then the set $\mathcal{B}=\{e_{ij}, f_{ij}: i,j=1,2,\ldots n\}$ constitutes a $B$-orthogonal basis of $\fr{u}(n)$.  Note that $e_{ij} = -e_{ji}$ and $f_{ij} = f_{ji}$.  Also, the Lie bracket relations among $e_{ij}$ and $f_{ij}$ are gives as follows:

\begin{lemma}\label{basebracketU}
The Lie brackets among the vectors {\rm (\ref{bracetSU})} satisfy the following relations:
\begin{eqnarray*}
[e_{ij}, e_{kl}] = \delta_{jk}e_{il} - \delta_{il}e_{kj} -\delta_{ik}e_{jl} -\delta_{jl}e_{ik}, &&
\left[ f_{ij}, f_{kl}\right] = -\delta_{jk}e_{il} + \delta_{il}e_{kj} -\delta_{ik}e_{jl} -\delta_{jl}e_{ik}, \nonumber\\
\left[ f_{ij}, e_{kl} \right] = \delta_{jk}f_{il} -\delta_{il}f_{kj} -\delta_{ik}f_{jl} -\delta_{jl}f_{ik}.
\end{eqnarray*}
\end{lemma} 

\noindent \underline{{\bf Case} {\bf (III)}} $G/H = \Sp(3)/\Sp(1)^3$.  We consider the embedding 
$$
\fr{h} = \fr{sp}(1)\oplus\fr{sp}(1)\oplus\fr{sp}(1) \hookrightarrow 
{\footnotesize \left\{ \begin{pmatrix}
ai  & 0 & 0 & \vline & -\bar{a}_1 & 0   & 0\\
0   & bi & 0 & \vline & 0  & -\bar{b}_1 & 0\\
0  & 0 & ci & \vline & 0  & 0 & -\bar{c}_1\\
\hline
a_1 & 0  & 0 & \vline & -ai &  0  & 0\\
0   & b_1 & 0 & \vline & 0 & -bi & 0\\
0  & 0 & c_1&  \vline & 0  & 0 & -ci
 \end{pmatrix}\, \Big|\, a,b,c\in\bb{R}, \ a_i, b_i, c_i\in\bb{C}  \right\} } \in\fr{sp}(3).
$$
Then, the tangent space $\fr{m} \cong T_{eH}(G/H)$ is given by $\fr{h}^{\perp}$ in $\fr{sp}(3)$ with respect to $B : \fr{sp}(3)\times\fr{sp}(3) \to \bb{R}$, $B(X, Y) = -\tr{XY}$.  It is easy to see that $\fr{m}$ has the following form:
$$
{\footnotesize \left\{ \begin{pmatrix}
0 & x_1 & x_2 & \vline & 0 &  -\bar{y}_1  & -\bar{y}_2 \\
-\bar{x}_1  & 0 & x_3 & \vline & -\bar{y}_1  & 0 & -\bar{y}_3\\
-\bar{x}_2  & -\bar{x}_3   & 0 & \vline & -\bar{y}_2  &-\bar{y}_3 & 0\\
\hline
0 & y_1  & y_2 & \vline & 0 &  \bar{x}_1  & \bar{x}_2\\
y_1   & 0 & y_3 & \vline & -x_1 & 0 & \bar{x}_3\\
y_2  & y_3 & 0 &  \vline & -x_2  & -x_3 & 0
 \end{pmatrix}\, \Big| \, x_i, y_i\in\bb{C}, i=1,2,3  \right\}. } 
$$
Let $M_{2n}\mathbb{C}$ be the set of $2n\times 2n$ complex matrices and we consider the following matrices in  $M_{2n}\mathbb{C}$:
\begin{itemize}
\item $E_{ab}$ with $1$ in $(a, b)$-entry and $-1$ in $(n+a, n+b)$-entry,
\item $F_{ab}$ with $i$ in $(a, b)$-entry and $i$ in $(n+a,$ $n+b)$-entry,
\item $G_{ab}$ with $-1$ in $(a, n+b)$-entry and $1$ in $(n+b, a)$-entry, 
\item $H_{ab}$ with $i$ in $(a, n+b)$-entry and $i$ in $(n+b, a)$-entry,
\end{itemize}
For $1\le a <b\le 2n$ we set
\begin{equation}\label{VectorsSP}
e_{ab}=E_{ab}-E_{ba}, \, f_{ab}=F_{ab}+F_{ba}, \, g_{ab}=G_{ab}+G_{ba}, \, h_{ab}=H_{ab}+H_{ba}.   
\end{equation}
Then the set $\mathcal{A}=\{e_{ab}, f_{ab}, g_{ab}, h_{ab} : 1\leq a<b\leq n, \, f_{aa}, g_{aa}, h_{aa} : 1\leq a\leq n\}$ constitutes a basis of $\fr{sp}(n)$, which is orthogonal with respect to $B$.  We have the following useful lemma   
\begin{lemma}\label{LemmaLieSP}
For the vectors {\rm (\ref{VectorsSP})} the following Lie bracket relations hold: 
\begin{center}
\begin{tabular}{lllll}
$[e_{ab},e_{bc}]=e_{ac}$  &  $[e_{ab},f_{bc}]=f_{ac}$ & $[e_{ab},g_{bc}]=g_{ac}$ & $[e_{ab},h_{bc}]=h_{ac}$ & $[f_{ab},f_{bc}]=-e_{ac}$ \\
$[f_{ab},g_{bc}]=-h_{ac}$ & $[f_{ab},h_{bc}]=g_{ac}$   &  $[g_{ab},g_{bc}]=-e_{ac}$   & $[g_{ab},h_{bc}]=-f_{ac}$  & $[h_{ab},h_{bc}]=-e_{ac}$ 
 \end{tabular}
\end{center}  
\end{lemma}

\subsection{Stiefel manifolds} 
The Stiefel manifold $V_{k}\bb{R}^n$ is the set of all orthonormal $k$-frames in $\bb{R}^n$.  This space is diffeomorphic to the reductive homogeneous space $\SO(n)/\SO(n-k)$.  We embed the group $\SO(n-k)$ in $\SO(n)$ as $\begin{pmatrix}
1_k & 0\\
0 & C
\end{pmatrix}$, where $C\in\SO(n-k)$.  The negative of the Killing form of $\fr{so}(n)$ is $B(X, Y)=-(n-2)\tr XY$.
Then, with respect to $B$, the subspace $\fr{m} =\fr{so}(n-k)^{\perp}$ in $\fr{so}(n)$ may be identified with the set  of matrices of the form
$$
\left\lbrace \begin{pmatrix}
D_k & A\\
-A^t & 0_{n-k}
\end{pmatrix} : D_k \in \fr{so}(k), A\in M_{k\times(n-k)}(\bb{R}) \right\rbrace.
$$
Let $E_{ab}$ denotes the $n\times n$ matrix with $1$ at the $(ab)$-entry and $0$ elsewhere.  Then the set $\mathcal{B}=\{\xi_{ab}=E_{ab}-E_{ba}: 1\le a\le k,\ 1\le a<b\le n\}$ constitutes a $B$-orthogonal basis of $\fr{m}$. 

\smallskip
In order to study $G$-invariant metrics on Stiefel manifolds $G/H$ we need to compute the isotropy representation.  Let $\varphi_n$ denotes the standard representation of $\SO(n)$ (given by the natural action of $\SO(n)$ on $\mathbb{R}^n$) and let $\wedge ^2\varphi _n$ be the second exterior power of $\varphi _n$.  Due to the isomorphism $\fr{so}(n)\cong \wedge^{2}\bb{R}^n$, it follows that $\Ad ^{\SO(n)}=\wedge ^2\varphi _n$.
The isotropy representation $\chi : \SO(n)\to \Aut(\fr{m})$ ($\fr{m}\cong T_{o}(G/H)$) of $G/H$ is characterized by the property
$\left.\Ad ^{\SO(n)}\right |_{\SO(n-k)}=\Ad ^{\SO(n-k)}\oplus\chi$. 
We compute
\begin{equation}\label{isotropy1}
\left.\Ad ^{\SO(n)}\right |_{\SO(n-k)}=\wedge ^2\varphi _n\big| _{\SO(n-k)}=\wedge ^2 (\varphi _{n-k}\oplus k)=\wedge ^2\varphi _{n-k}\oplus \wedge ^2\mathbbm{1}_k\oplus (\varphi _{n-k}\wedge \mathbbm{1}_k),
\end{equation}
where $\wedge^{2}\mathbbm{1}_k$ is the sum of ${k}\choose{2}$ trivial representations $\mathbbm{1}$.  Therefore, the isotropy representation is given by $\chi = \mathbbm{1}\oplus \cdots\oplus \mathbbm{1}\oplus\varphi _{n-k}\oplus \cdots\oplus\varphi _{n-k}$.  This decomposition induces an $\Ad^{H}$-invariant decomposition of $\fr{m}$ given by $\fr{m}=\fr{m}_1\oplus\cdots\oplus\fr{m}_s$,
where the first ${k}\choose{2}$ $\Ad^{H}$-modules are $1$-dimensional and the remaining $k$ modules are $(n-k)$-dimensional.  It is clear that the isotropy representation of $V_k\bb{R}^n$ contains equivalent summands, so a complete description of all $G$-invariant metrics is rather hard.  In \cite{ADN} the authors introduced a method for proving existence of homogeneous Einstein metrics by assuming additional symmetries.  In \cite{Sta1} the author presented a systematic and organized description of such metrics.

\subsubsection{The Stiefel manifold $V_2\bb{R}^{n}\cong \SO(n)/\SO(n-2)$}
The isotropy representation of $V_2\bb{R}^n$ is expressed as a direct sum $\chi = \mathbbm{1}\oplus\chi_1\oplus\chi_2$, where $\chi_1\approx \chi_2= \varphi_{n-2}$ is the standard representation of $\SO(n-2)$.  This decomposition induces an $\Ad^{\SO(n-2)}$-invariant decomposition of $\fr{m}$ given by $\fr{m} = \fr{m}_{0}\oplus\fr{m}_1\oplus\fr{m}_2$.  
%
%
It can be shown (cf. \cite{Ke}) that any $\SO(n)$-invariant metric on $V_2\bb{R}^n$ (modulo isometries) can be described by an $\Ad^{\SO(n-2)}$-invariant inner product on $\fr{m}$ of the form 
$\Lambda= \lambda_0{\rm Id}_{\fr{m}_0} + \lambda_1{\rm Id}_{\fr{m}_1} +\lambda_2{\rm Id}_{\fr{m}_2},$ $\lambda_i\in\bb{R}_{+}, i =0,1,2$.

\begin{remark}
\textnormal{Using Lemma \ref{bracSO} it is easy to see that modules $\fr{m}_{i}$, $i = 0,1,2$ satisfy the Lie bracket relation $[\fr{m}_{i}, \fr{m}_{i}]\subset\fr{so}(n-2)$.  Therefore, the Stiefel manifold $V_2\bb{R}^n$ is also a generalized Wallach space.}
\end{remark}

\begin{theorem}\textnormal{(\cite{Arv2}, \cite{Ke})}\label{theoremEinstein}
The Stiefel manifold $V_2\bb{R}^n = \SO(n)/\SO(n-2)$ admits (up to scale) exactly one $\SO(n)$-invariant Einstein metric which is given explicitly as $(1, \lambda_1, \lambda_1)$ where $\lambda_1 = ({n-1})/{2(n-2)}$.
\end{theorem}

\subsubsection{The Stiefel manifolds $V_{1+n_2}\bb{R}^n$}
Let $G/H$ be the Stiefel manifold $V_{1+n_2}\bb{R}^n = \SO(n)/\SO(n_3)$, with $n = 1+n_2+n_3$.  The isotropy representation on this case according to (\ref{isotropy1}) contains equivalent summands.  We will describe a special class of invariant metrics in this space (for more details see for example \cite{Sta1}, \cite{ArSaSt1} and \cite{ArSaSt2}).  The basic approach is to use an appropriate subgroup $K$ of $G$, such that the special class of $\Ad^{K}$-invariant inner products, which are a subset of $\Ad^{H}$-invariant inner products, are diagonal.  In order to have this, it is sufficient for the subgroup $K$ to satisfy the condition $H\subset K \subset N_{G}(H)\subset G$.

We take the subgroup $K = \SO(n_2)\cdot\SO(n_3)$ of $\SO(n)$.  Then, for the tangent space $\fr{m} \cong T_{o}(G/H)$, we consider the irreducible, $\Ad^{K}$-invariant and non-equivalent decomposition: 
$
\fr{m} = \fr{so}(n_2)\oplus  \fr{m}_{12} \oplus  \fr{m}_{13} \oplus  \fr{m}_{23}. 
$
Therefore, the $G$-invariant metrics on $G/H$ determined by the $\Ad^{K}$-invariant inner products on $\fr{m}$ which are are given by 
$
\Lambda = \lambda_2 {\rm Id}_{ \fr{so}(n_2)} +\lambda_{12} {\rm Id}_{ \fr{m}_{12}} + \lambda_{13} {\rm Id}_{ \fr{m}_{13}} + \lambda_{23} {\rm Id}_{ \fr{m}_{23}}
$, $\lambda_2, \lambda_{ij}\in\bb{R}_{+}$, $1\leq i < j \leq 3$.

\begin{remark}\label{remarkJensen}
\textnormal{
It was proved in \cite{ArSaSt1} that the Stiefel manifolds $V_{4}\bb{R}^n$ $(n=1+n_2+n_3\geq 6)$ for $n_2 = 3, n_3 = n-4$, admit at least four invariant Einstein metrics.  Two of them are Jensen's Einstein metrics, which have the form $\Lambda_{\tiny\mathsf{Jen}} = (\lambda_2, \lambda_2, 1, 1),$ where $\lambda_2 = \displaystyle{\frac{-2+n \pm \sqrt{7-7n+n^2}}{n-1}}$. 
}
\end{remark}

\subsection{Spheres}
Spheres can be expressed as homogeneous spaces in various different ways.  More precisely, we have the isotropy irreducible spheres $\bb{S}^{n} \cong \SO(n+1)/\SO(n), \bb{S}^6 \cong \G_2/\SU(3), \bb{S}^{7}\cong {\rm Spin(7)}/\G_2$, the spheres for which the isotropy representation splits into two irreducible and non equivalent modules $\bb{S}^{2n+1}\cong \U(n+1)/\U(n) \cong \SU(n+1)/\SU(n),$ $\bb{S}^{4n+1} \cong \Sp(n+1)/\Sp(n) \cong (\Sp(n+1)\cdot\Sp(1))/\Sp(n)\cdot\Sp(1)$, $\bb{S}^{15} \cong {\rm Spin(9)}/{\rm Spin(7)}$ and  the sphere $\bb{S}^{4n+1} \cong (\Sp(n+1)\cdot\U(1))/\Sp(n)\cdot \U(1)$ with three irreducible and non equivalent isotropy modules. 

\subsubsection{The sphere $\U(n+1)/\U(n)\cong \bb{S}^{2n+1}$}\label{sectionSphere1}
We embed the subgroup $\U(n)$ into $\U(n+1)$ as 
$
\begin{pmatrix}
1 & 0\\
0 & A 
\end{pmatrix},
$
where $A\in\U(n)$.  We consider the $\Ad^{\U(n+1)}$-invariant inner product $B : \fr{u}(n+1)\times\fr{u}(n+1) \to \bb{R}$ given by $B(X, Y)=-\tr XY$.  Then, with respect to this inner product, the subspace $\fr{m} =\fr{u}(n)^{\perp}$ in $\fr{u}(n+1)$ may be identified with the set  of matrices of the form
$$
{\small \left\lbrace \begin{pmatrix}
i \al  & a_1 & \ldots & a_{n}\\
-\bar{a}_1  & 0 &  &  0 \\
\vdots & \vdots &  \ddots &  \vdots\\
-\bar{a}_{n} &0 & \ldots & 0
\end{pmatrix} : \al\in \bb{R}, a_i \in \bb{C}, i=1,2,\ldots, n  \right\rbrace. }
$$ 
Let $E_{ab}$ denote the $n\times n$ matrix with $1$ at the $(ab)$-entry and $0$ elsewhere.  We set
$
e_{1j} = E_{1j}-E_{j1}, \ f_{1j} = i(E_{1j} + E_{j1}).
$
Then the set $\mathcal{B}=\{e_{1j}, f_{1j}: j=1,2,\ldots n+1\}$ constitutes a $B$-orthogonal basis of $\fr{m}$. 

We compute the isotropy representation $\chi\otimes\bb{C} : \U(n) \to \Aut(\fr{m}\otimes\bb{C})$.   
 If is $\mu_{n}$ the standard representation of $\U(n)$, then we have that $\Ad^{\U(n)}\otimes\bb{C} = \mu_{n}\otimes\bar{\mu}_{n}$.  From Proposition \ref{isotrepr} it is
\begin{eqnarray*}
\Ad^{\U(n+1)}\otimes\bb{C}\big|_{\U(n)} &=& (\mu_{n+1}\otimes\bar{\mu}_{n+1})\big|_{\U(n)} = (\mu_{n}\oplus\mathbbm{1})\otimes\overline{{(\mu_{n}\otimes\mathbbm{1})}} = (\mu_{n}\otimes\bar{\mu}_{n}) \oplus \mu_{n}\oplus\bar{\mu}_{n}\oplus(\mathbbm{1}\otimes\mathbbm{1})\\
&=& \Ad^{\U(n)}\otimes\bb{C} \oplus \mu_{n}\oplus\bar{\mu}_{n}\oplus \mathbbm{1},
\end{eqnarray*}
therefore $\chi\otimes\bb{C} = \mathbbm{1}\oplus\mu_{n}\oplus\bar{\mu}_{n}$.  This decomposition induces an $\Ad^{H}\otimes\bb{C}$-invariant decomposition of $\fr{m}\otimes\bb{C}$ into three irreducible $\Ad^{H}\otimes\bb{C}$-modules, that is $\fr{m}\otimes\bb{C}=\fr{p}_1\oplus\fr{p}_2\oplus\fr{p}_3$, where $\dim_{\bb{C}}\fr{p}_i = n$, $(i=2,3)$ and $\dim_{\bb{C}}\fr{p}_1 = 1$.  For the real tangent space $\fr{m}$ we have $\fr{m} = \fr{m}_1\oplus\fr{m}_2$, where $\fr{m}_{1}\otimes\bb{C} = \fr{p}_1$, $\fr{m}_2\otimes\bb{C} = \fr{p}_2\oplus\fr{p}_3$ are non equivalent and the dimensions are $\dim\fr{m}_1 = 1$ and $\dim\fr{m}_2 = 2n$.
The $\U(n+1)$-invariant metrics are determined by the endomorphism $\Lambda : \fr{m} \to \fr{m}$ which has diagonal form
$
\Lambda = \lambda_1{\rm Id}_{\fr{m}_1} + \lambda_2{\rm Id}_{\fr{m}_2}, \  \lambda_{i}\in\bb{R}^{+}, i=1,2.
$

\subsubsection{The sphere $\Sp(n+1)/\Sp(n)\cong \bb{S}^{4n+3}$}\label{sectionSphere2}
We embed the Lie subalgebra \begin{equation*}
\mathfrak{sp}(n)=\left\{\begin{pmatrix}
X & -{}\bar{Y}\\
Y & \bar{X}
\end{pmatrix}\ \Big\vert 
 \begin{array}{l}X\in\mathfrak{u}(n), \\ 
Y\ \mbox{ is \ a}\  n\times n\ \mbox{complex symmetric matrix}
\end{array}
\right\}\ 
\end{equation*} 
 in the Lie algebra $\mathfrak{sp}(n+1)$ as follows: 
 {\footnotesize
$
  \left\{\begin{pmatrix}
 0  & 0 &\vline & 0   & 0\\
0  & X & \vline & 0  & -\bar{Y}\\
\hline
0  & 0 &\vline & 0   & 0\\
0  & Y & \vline & 0  & \bar{X}
 \end{pmatrix} \right\}. 
 $ }
We consider the $\Ad^{\Sp(n+1)}$-invariant inner product $B : \fr{sp}(n+1)\times\fr{sp}(n+1) \to \bb{R}$ given by $B(X, Y)=-\tr XY$.  Then, with respect to this inner product, the subspace $\fr{m} =\fr{sp}(n)^{\perp}$ in $\fr{sp}(n+1)$, may be identified with the set  of matrices of the form

{\footnotesize $$
\left\lbrace \begin{pmatrix}
ia  & \al_1 & \cdots & \al_n &\vline & -\bar{b}_1 & -\bar{b}_2 & \cdots & -\bar{b}_{n+1} \\
-\bar{\al}_1 & 0 & \cdots  & 0 & \vline & -\bar{b}_2 & 0 & \cdots & 0\\
\vdots & \vdots &\ddots   & \vdots & \vline & \vdots & \vdots &  \ddots &  \vdots\\
-\bar{\al}_n  & 0 & \cdots & 0 & \vline & -\bar{b}_{n+1} &  0 & \cdots & 0\\
\hline
b_{1}  & b_2 & \cdots & b_{n+1} &\vline & -ia & \bar{\al}_1 & \cdots & \bar{\al}_{n} \\
b_2 & 0 & \cdots  & 0 & \vline & -\al_1 & 0 & \cdots & 0\\
\vdots & \vdots &\ddots   & \vdots & \vline & \vdots & \vdots &  \ddots &  \vdots\\
b_{n+1}  & 0 & \cdots & 0 & \vline & -\al_n &  0 & \cdots & 0\\
 \end{pmatrix} : 
 \begin{array}{l} a\in\bb{R}, \, \al_i\in\bb{C}, i=1,2,\ldots n \\ 
b_j \in\bb{C}, j = 1,2,\ldots, n+1
\end{array}
\right\rbrace.
$$ }
Let $M_{2n}\mathbb{C}$ be the set of $2n\times 2n$ complex matrices and we consider the following matrices in  $M_{2n}\mathbb{C}$:
\begin{itemize}
\item $E_{1b}$ with $1$ in $(1, b)$-entry and $-1$ in $(n+2, n+1+b)$-entry,
\item $F_{1b}$ with $i$ in $(1, b)$-entry and $i$ in $(n+2,$ $n+1+b)$-entry,
\item $G_{1b}$ with $-1$ in $(1, n+1+b)$-entry and $1$ in $(n+1+b, 1)$-entry, 
\item $H_{1b}$ with $i$ in $(1, n+1+b)$-entry and $i$ in $(n+1+b, 1)$-entry,
\item $u_{11}$ with $i$ in $(1, 1)$-entry and $-i$ in $(n+2, n+2)$-entry,
\item $k_{1n+2}$ with $-1$ in $(1, n+2)$-entry and $1$ in $(n+2, 1)$-entry, 
\item $\mu_{1n+2}$ with $i$ in $(1, n+2)$-entry and $i$ in $(n+2, 1)$-entry.
\end{itemize}
For $1\le a <b\le 2n$ we set
$
e_{1b}=E_{1b}-E_{b1}, \, f_{1b}=F_{1b}+F_{b1}, \, g_{1b}=G_{1b}+G_{b1}, \,
h_{1b}=H_{1b}+H_{b1}.   
$ 
Then the set $\mathcal{A}=\{u_{11}, k_{1n+2}, \mu_{1n+2}, e_{12}, f_{12},\ldots, e_{1n+1}, f_{1n+1}, g_{12}, h_{12},\ldots, g_{1n+1}, h_{1n+1}\}$ constitutes a basis of $\fr{m}$, which is orthogonal with respect to $B$.  For the Lie brackets among the previous vectors we have  Lemma \ref{LemmaLieSP} together with  the following:

\begin{lemma}\label{LemmaLieSP2}
For the vectors $u_{11}, k_{1n+2}, \mu_{1n+2}, e_{ab}, f_{ab}, g_{ab}, h_{ab}$ the following Lie bracket relations hold: 
\begin{tabular}{lllll}
$[u_{11}, e_{ab}] = f_{ab}$ & $[u_{11}, f_{ab}] = -e_{ab}$ & $[u_{11}, k_{1n+2}] = -2\mu_{1n+2}$ & $[u_{11}, \mu_{1n+2}] = 2k_{1n+2}$ & $[u_{11}, g_{ab}] = -h_{ab}$\\
$[u_{11}, h_{ab}] = g_{ab}$ & $[e_{ab}, k_{1n+2}] = -g_{ab}$ & $[e_{ab}, \mu_{1n+2}] = -h_{ab}$ & $[f_{ab}, k_{1n+2}] = -h_{ab}$ & $[f_{ab}, \mu_{1n+2}] = g_{ab}$ \\
$[g_{ab}, k_{1n+2}] = e_{ab}$ & $[g_{ab}, \mu_{1n+2}] = -f_{ab}$ & $[h_{ab}, k_{1n+2}] = f_{ab}$ & $[h_{ab}, \mu_{1n+2}] = e_{ab}.$  &  
 \end{tabular}
\end{lemma}

We compute the isotropy representation $\chi\otimes\bb{C} : \Sp(n) \to \Aut(\fr{m}\otimes\bb{C})$.  
If $\nu_{n}$ is the standard representation of $\Sp(n)$, then we have that $\Ad^{\Sp(n)}\otimes\bb{C} = S^2\nu_{n}$.  From Proposition \ref{isotrepr} it is
\begin{eqnarray*}
\Ad^{\Sp(n+1)}\otimes\bb{C}\big|_{\Sp(n)} &=& S^2\nu_{n+1}\big|_{\Sp(n)} = S^{2}(\nu_{n}\oplus  \mathbbm{1}) = S^2\nu_n\oplus S^{2} \mathbbm{1}\oplus(\nu_n\otimes  \mathbbm{1})\\
&=& \Ad^{\Sp(n)}\otimes\bb{C} \oplus S^{2} \mathbbm{1}\oplus(\nu_n\otimes  \mathbbm{1}),
\end{eqnarray*}

therefore $\chi \otimes\bb{C}= S^2\mathbbm{1}\oplus\nu_{n}$.  This decomposition induces an $\Ad^{\Sp(n)}\otimes\bb{C}$-invariant decomposition of $\fr{m}\otimes\bb{C}$ into two irreducible $\Ad^{\Sp(n)}\otimes\bb{C}$-modules, that is $\fr{m}\otimes\bb{C}=\fr{p}_1\oplus\fr{p}_2$.  For the real tangent space $\fr{m}$ we have $\fr{m} = \fr{m}_1\oplus\fr{m}_2$, where $\fr{m}_{1}\otimes\bb{C} = \fr{p}_1$, $\fr{m}_2\otimes\bb{C} = \fr{p}_2$ are non equivalent and the dimensions are $\dim\fr{m}_1 = 3$ and $\dim\fr{m}_2 = 4n$.
The $\Sp(n+1)$-invariant metrics are given by the endomorphism $\Lambda : \fr{m} \to \fr{m}$, which has diagonal form
$
\Lambda = \lambda_1{\rm Id}_{\fr{m}_1} + \lambda_2{\rm Id}_{\fr{m}_2}, \  \lambda_{i}\in\bb{R}^{+}, i=1,2.
$

\section{Equigeodesics}\label{section3}
Let $(G/H, g)$ be a Riemannian homogeneous space with $G$ compact and semisimple and let $\Lambda : \fr{m}\to \fr{m}$ be the $\Ad^H$-equivariant, positive definite, symmetric operator associated to $g$, defined by $\langle\cdot, \cdot\rangle = B(\Lambda\cdot, \cdot)$  (cf. subsection \ref{subsection1.3}).  In the following we will use the same notation for $g$ and $\Lambda$.
\begin{definition}\label{def1}
 A geodesic $\gamma(t)$ on $(G/H, g)$ through the origin $o=eH$ is called homogeneous if it is the orbit of a one-parameter subgroup of $G$, that is $\gamma(t) = (\exp tX )\cdot o$, where $X\in\fr{g}$.  The vector $X$ is called a geodesic vector.
\end{definition}
Definition \ref{def1} establishes a one-to-one correspondence between geodesic vectors $X$ and homogeneous geodesics at the origin.  It is well known that every homogeneous Riemannian manifold admits at least one homogeneous geodesic \cite{KoSz}. In the case where the group $G$ is semisimple we have the following existence result:
\begin{theorem}{\rm (\cite{KoSz})}
If $G$ is semisimple then a Riemannian manifold $G/H$ admits at least $n = \dim(G/H)$ mutually orthogonal homogeneous geodesics through the origin $o=eH$.
\end{theorem}

The next theorem (due to Kowalski and Vanhecke \cite{KoVa}) gives an algebraic characterization of the geodesic vectors.
\begin{theorem}
If $g$ is a $G$-invariant metric, then a non zero vector $X\in\fr{g}$ is a geodesic vector if and only if
$\langle X_{\fr{m}}, [X,Z]_\fr{m}\rangle = 0,$ for all $Z\in\fr{m}$.
\end{theorem}

A curve of the form $\gamma(t) = (\exp tX)\cdot o$ is called {\it equigeodesic} on $G/H$ if it is a geodesic with respect to each invariant metric $\Lambda$.  The vector $X$ is called {\it equigeodesic vector}.  The following proposition gives  an algebraic characterization of equigeodesic in terms of equigeodesic vectors, that is vectors $X\in\fr{m}$ such that the orbit $\gamma(t) = (\exp tX)\cdot o$ is an homogeneous equigeodesic.

\begin{prop}\textnormal{(\cite{CGN})}\label{sxesi}
Let $G/H$ be a reductive homogeneous space with reductive decomposition $\fr{g}=\fr{h}\oplus\fr{m}$ and $X\in\fr{m}$ be a non-zero vector.  Then $X$ is an equigeodesic vector if and only if
\begin{equation}\label{equigeod}
[X, \Lambda X]_{\fr{m}}=0,
\end{equation}
for each invariant metric $\Lambda$.
\end{prop} 

The process of solving equation (\ref{equigeod}) is equivalent to solving a non linear algebraic system of equations whose variables are the coefficients of the vector $X$.   In some cases (depending just on the $\fr{m}_i$-parts of $X$) the non linear system vanishes completely (i.e. the system is identically zero).  This motivates the following definition.

\begin{definition}
An equigeodesic vector $X$ is said to be

\smallskip
$(a)$ structural if the algebraic system associated to equation {\rm (\ref{equigeod})} is identically zero.

\smallskip
$(b)$ algebraic if the coordinates of the vector $X$ come from a solution of a (not identically zero) non linear algebraic system associated to equation {\rm (\ref{equigeod})}. 
\end{definition}

{From the invariance of the metric $\Lambda$, we have that $\Lambda|_{\fr{m}_i} = \lambda_{i}{\rm Id}_{\fr{m}_i}$, for some $\lambda_i > 0$, for each irreducible component of the isotropy representation.  Therefore, if $X\in\fr{m}_i$ then equation (\ref{equigeod}) is satisfied trivially.  }

\begin{definition}
An equigeodesic vector $X\in\fr{m}$ called

\smallskip
$(a)$ trivial if $X\in\fr{m}_i$ for some $i$,

\smallskip
$(b)$ non trivial if $X\in\bigoplus_{i=1}^s\fr{m}_{i}$ with $s>1$  and $X$ has components in at least two of the $\fr{m}_i$'s. 
\end{definition}

\begin{remark}
{\rm It is obvious that the trivial equigeodesic vectors are structural equigeodesic vectors.} 
\end{remark}

\begin{lemma}\label{Lemma3}
Let $G/H$ be a homogeneous space admitting a reductive decomposition $\fr{g} = \fr{h}\oplus\fr{m}$.  Assume that $\fr{m}$ decomposes into a direct sum 
$$
\fr{m} = \bigoplus_{i=1}^{s}\fr{m}_{i}
$$
of non equivalent $\Ad^{H}$-invariant summands.  Then a vector $X=\sum_{i=1}^{s}X_{\fr{m}_i}$, with $X_{\fr{m}_{i}}\in\fr{m}_i$ for $1\leq i\leq s$ is equigeodesic if and only if
$$
[X_{\fr{m}_i}, X_{\fr{m}_j}]_{\fr{m}} = 0, \ \ \mbox{for} \ \ 1\leq i< j\leq s.
$$
\end{lemma}
\begin{proof}
If $\pi : \fr{g} \to \fr{m}$ is the projection onto $\fr{m}$, then $\pi([X, \Lambda X])=[X, \Lambda X]_{\fr{m}}$.  Let $X=\sum_{i=1}^{s}X_{\fr{m}_i}\in\fr{m}$.  Then 
\begin{eqnarray*}
&&[X, \Lambda X]_{\fr{m}} = \pi([X, \Lambda X]) = \pi(\big[\sum_{i=1}^{s}X_{\fr{m}_i}, \Lambda(\sum_{i=1}^{s}X_{\fr{m}_{i}})\big]) = \pi(\big[\sum_{i=1}^{s}X_{\fr{m}_{i}}, \sum_{i=1}^{s}\lambda_{i}X_{\fr{m}_{i}}\big])\\
&& = \sum_{i=2}^{s}(\lambda_{i}-\lambda_1)\pi([X_{\fr{m}_1}, X_{\fr{m}_i}]) + \sum_{i=3}^{s}(\lambda_i-\lambda_2)\pi([X_{\fr{m}_2}, X_{\fr{m}_i}]) + \cdots + (\lambda_s - \lambda_{s-1})\pi([X_{\fr{m}_{s-1}}, X_{\fr{m}_s}]) \\
&& = \sum_{i=2}^{s}(\lambda_{i}-\lambda_1)[X_{\fr{m}_1}, X_{\fr{m}_i}]_{\fr{m}} + \sum_{i=3}^{s}(\lambda_i-\lambda_2)[X_{\fr{m}_2}, X_{\fr{m}_i}]_{\fr{m}} + \cdots + (\lambda_s - \lambda_{s-1})[X_{\fr{m}_{s-1}}, X_{\fr{m}_s}]_{\fr{m}}.
\end{eqnarray*}
We know that $X$ is an equigeodesic vector if and only if $[X, \Lambda X]_{\fr{m}} = 0$ for each invariant metric $\Lambda = (\lambda_1,\ldots,\lambda_s)$ $\lambda_i\in\bb{R}^{+}, i=1,2,\ldots,s$.  This occurs if and only if $[X_{\fr{m}_i}, X_{\fr{m}_j}]_{\fr{m}} = 0$, where $1\leq i<j \leq s$.  
\end{proof}

Some classes of structural equigeodesic vectors were given for generalized flag manifolds with two and $s>2$ isotropy summands (\cite{GrNe}, \cite{WaZh}).  Also, in \cite{Sta2} the author described all structural equigeodesic vectors for the generalized flag manifolds with $G_2$-type $\fr{t}$-roots, namely $\F_4/\U(3)\cdot\U(1)$, $\E_6/\U(3)\cdot\U(3)$, $\E_7/\U(6)\cdot\U(1)$ and gave partial description for $\E_8/\E_6\cdot\U(1)\cdot\U(1)$. 


\section{Algebraic Equigeodesic Vectors on Reductive Homogeneous Spaces}

In this section we describe  algebraic equigeodesic vectors for certain reductive homogeneous spaces, namely $V_2\bb{R}^n$, $W^{6} = \U(3)/\U(1)^3$, $W^{12} = \Sp(3)/\Sp(1)^3$ and spheres $\bb{S}^{2n+1} = \U(n+1)/\U(n)$, $\bb{S}^{4n+3} = \Sp(n+1)/\Sp(n)$.  Also, for the spaces $\SO(n_1+n_2+n_3)/\SO(n_1)\cdot\SO(n_2)\cdot\SO(n_3)$, $V_{1+n_2}\bb{R}^{n}$ and $\bb{S}^{4n+3} = \Sp(n+1)/\Sp(n)$ we obtain the general algebraic non linear system of equations whose variables are the coefficients of the vector.  Since these systems are more complicated, we examine in full detail only some examples of such spaces, such as $\SO(6)/\SO(3)\cdot\SO(2)$ and for $V_4\bb{R}^6$ with respect to Jensen's Einstein metrics.  


\subsection{Algebraic equigeodesic vectors on generalized Wallach spaces} 

\subsubsection{The genaralized Wallach space $\SO(n)/\SO(n_1)\cdot\SO(n_2)\cdot\SO(n_3)$}
We consider the generalized Wallach space $G/H = \SO(n)/\SO(n_1)\cdot\SO(n_2)\cdot\SO(n_3)$ with $n = n_1+n_2+n_3$.  Then the tangent space $T_{o}(G/H)$ is written as $\fr{m} = \bigoplus_{1\leq i<j \leq 3}\fr{m}_{ij}$ and for each module we have the following description:
\begin{eqnarray}\label{baseWallach}
\fr{m}_{ij} = \Span\{\xi_{ab} : n_1 +\cdots +n_{i-1} +1 \leq a \leq n_1 + \cdots + n_i, 
 n_1+\cdots + n_{j-1} < b\leq n_1 + \cdots + n_{j} \}. &&
\end{eqnarray}

The $\SO(n)$-invariant metrics $\Lambda_{{\tiny\mathsf{Wal}}}$ on $G/H$ have a diagonal form and are given as $\Lambda_{{\tiny\mathsf{Wal}}}|_{\fr{m}_{ij}} = \lambda_{ij}{\rm Id}_{\fr{m}_{ij}}$.  For $X_{\fr{m}_{ij}}\in\fr{m}_{ij}$  we write
\begin{eqnarray}\label{vectorW}
X_{\fr{m}_{ij}} = \sum_{k\in I_i, l\in J_j} a_{kl} \xi_{kl}, \ \ a_{kl}\in\bb{R},
\end{eqnarray}
where $I_{i} = \{n_1 +\cdots +n_{i-1} +1 \leq k \leq n_1 + \cdots + n_i\}, 
J_{j} =\{ n_1+\cdots + n_{j-1} < l \leq n_1 + \cdots + n_{j} \} $ 
and $a_{kl}\in\bb{R}$.  We have

\begin{prop}\label{PropWallach}
Let the generalized Wallach space $\SO(n)/\SO(n_1)\cdot\SO(n_2)\cdot\SO(n_3)$ and let $X = X_{\fr{m}_{12}} + X_{\fr{m}_{13}} + X_{\fr{m}_{23}}\in\fr{m}$, where $X_{\fr{m}_{ij}}$ are as in {\rm (\ref{vectorW})}.  Then the vector $X$ is an euigeodesic vector if and only if  
\begin{eqnarray}\label{systemW}
{\mathsf{(A)}}\left \{ \begin{array}{l} 
\sum_{i\in I_{1}, j\in J_{2}}\sum_{k\in I_{1}, l\in J_{3}} a_{ij}a_{kl} B([\xi_{ij}, \xi_{ab}]_{\fr{m}}, \xi_{kl})=0\\ 
 \quad \quad \quad \quad \quad \quad \quad \quad  \    \vdots\\ 
\sum_{i\in I_{1}, j\in J_{2}}\sum_{k\in I_{1}, l\in J_{3}} a_{ij}a_{kl} B([\xi_{ij}, \xi_{d_3}]_{\fr{m}}, \xi_{kl})=0\nonumber\\ 
\end{array} \right. \\
{\mathsf{(B)}}\left \{ \begin{array}{l} 
\sum_{i\in I_{1}, j\in J_2}\sum_{k\in I_2, l\in J_3} a_{ij}a_{kl} B([\xi_{ij}, \xi_{ab}]_{\fr{m}}, \xi_{kl})=0\\ 
 \quad \quad \quad \quad \quad \quad \quad \quad  \    \vdots \\ 
\sum_{i\in I_1, j\in J_2}\sum_{k\in I_2, l\in J_3} a_{ij}a_{kl} B([\xi_{ij}, \xi_{d_2}]_{\fr{m}}, \xi_{kl})=0\\ 
\end{array} \right. \\
{\mathsf{(C)}}\left \{ \begin{array}{l} 
\sum_{i\in I_1, j\in J_3}\sum_{k\in I_2, l\in J_3} a_{ij}a_{kl} B([\xi_{ij}, \xi_{ab}]_{\fr{m}}, \xi_{kl})=0\\ 
 \quad \quad \quad \quad \quad \quad \quad \quad \      \vdots\\ 
\sum_{i\in I_1, j\in J_3}\sum_{k\in I_2, l\in J_3} a_{ij}a_{kl} B([\xi_{ij}, \xi_{d_1}]_{\fr{m}}, \xi_{kl})=0.\nonumber 
\end{array} \right.
\end{eqnarray}



\end{prop}
\begin{proof}
Let $X, Y\in\fr{m}$, where $X = X_{\fr{m}_{12}}+X_{\fr{m}_{13}}+X_{\fr{m}_{23}}$ and let $B$ be the negative of the Killing form of $\SO(n)$.  We have
\begin{eqnarray}\label{proofW}
B([X, \Lambda_{{\tiny\mathsf{Wal}}} X]_{\fr{m}}, Y) &=& B([\sum_{1\leq i< j\leq 3}X_{\fr{m}_{ij}},\ \Lambda_{{\tiny\mathsf{Wal}}} \sum_{1\leq i< j\leq 3}X_{\fr{m}_{ij}}]_{\fr{m}},\, Y)  \nonumber\\
&=& B([\sum_{1\leq i< j\leq 3}X_{\fr{m}_{ij}}, \sum_{1\leq i< j\leq 3}\lambda_{ij}X_{\fr{m}_{ij}}]_{\fr{m}}, \, Y)
= (\lambda_{13} - \lambda_{12}) B([X_{\fr{m}_{13}}, X_{\fr{m}_{12}}]_{\fr{m}}, \, Y) \nonumber\\
&& + (\lambda_{23}-\lambda_{12}) B([X_{\fr{m}_{12}}, X_{\fr{m}_{23}}]_{\fr{m}}, \, Y)  + (\lambda_{23} - \lambda_{13}) B([X_{\fr{m}_{13}}, X_{\fr{m}_{23}}]_{\fr{m}}, \, Y).
\end{eqnarray}
Then $X$ is an equigeodesic vector if and only if for any $\SO(n)$-invariant metric (\ref{Nikonorov}) and for all $Y\in\fr{m}$ it is $B([X, \Lambda_{{\tiny\mathsf{Wal}}} X], Y) = 0$, that is $B([X_{\fr{m}_{13}}, X_{\fr{m}_{12}}]_{\fr{m}},\, Y) = B([X_{\fr{m}_{12}}, X_{\fr{m}_{23}}]_{\fr{m}}, \, Y) = B([X_{\fr{m}_{13}}, X_{\fr{m}_{23}}]_{\fr{m}}, \, Y) = 0$. 
For example, if we take $X_{\fr{m}_{12}} = \sum_{i\in I_1, j\in J_2} a_{ij} \xi_{ij}$ and $X_{\fr{m}_{13}} = \sum_{k\in I_1, l\in J_3} a_{kl} \xi_{kl}$ (as in (\ref{vectorW})), then from Lemma \ref{BrackW} it is $[X_{\fr{m}_{13}}, X_{\fr{m}_{12}}] \in\fr{m}_{23}$, therefore for $Y\in \fr{m}_{13}$ or $Y\in\fr{m}_{12}$ it is $B([X_{\fr{m}_{13}}, X_{\fr{m}_{12}}], Y) = 0$ (because the summands $\fr{m}_{ij}$ are pairwise orthogonal).  In the case where $Y\in\fr{m}_{23}$,  from (\ref{proofW}) we have that the second and third term are equal to zero, so we obtain  that
\begin{eqnarray*}
0=B([X_{\fr{m}_{13}}, X_{\fr{m}_{12}}]_{\fr{m}},\, Y) &=& B([X_{\fr{m}_{12}}, Y]_{\fr{m}}, \, X_{\fr{m}_{13}}) =  B([\sum_{i\in I_1,j\in J_2} a_{ij} \xi_{ij}, \ Y]_{\fr{m}}, \  \sum_{k\in I_1, l\in J_3} a_{kl} \xi_{kl})\\ 
&=& \sum_{i\in I_1,j\in J_2} \sum_{k\in I_1, l\in J_3}  a_{ij} a_{kl} B([\xi_{ij}, Y]_{\fr{m}}, \, \xi_{kl}).
\end{eqnarray*}
By substituting
$ 
Y= \xi_{ab}\in \fr{m}_{23} = \Span\{\xi_{ab} : n_1 +1 \leq a \leq n_1 + n_2,
 n_1 + n_{2} < b\leq n_1 + n_{2} + n_3 \},
$
into the above equation, we obtain system  {$\mathsf{(A)}$}.  By  similar calculations we obtain systems {$\mathsf{(B)}$} and  {$\mathsf{(C)}$}. 
\end{proof}

By using Lemma \ref{bracSO} we can compute the brackets relations in system (\ref{systemW}).  
Next, we will study the equigeodesic vectors on the generalized Wallach space $\SO(6)/\SO(3)\cdot\SO(2)$, that is when $n_1 = 1, n_2 = 3$ and $n_3 = 2$.   
Then system (\ref{systemW}) becomes
\begin{eqnarray}
\left \{ \begin{array}{l}
 a_{15}a_{12} = a_{15}a_{13} = a_{15}a_{14} = a_{16}a_{12} = a_{16}a_{13} = a_{16}a_{14} = 0\nonumber\\
 a_{25}a_{12}+a_{35}a_{13}+a_{45}a_{14} = 0, \ \ a_{26}a_{12} + a_{36}a_{13} + a_{46}a_{14} = 0\\
 a_{25}a_{15} + a_{26}a_{16} = 0,\ \ a_{35}a_{15} + a_{36}a_{16} = 0,\ \ a_{45}a_{15} + a_{46}a_{16}=0.\nonumber
\end{array} \right.
\end{eqnarray}
We take the following solutions:
\begin{eqnarray*}
&&\{a_{15} = 0, a_{16} = 0, a_{45}  -(a_{12} a_{25} + a_{13} a_{35})/a_{14}, a_{46} = -(a_{12} a_{26} + a_{13} a_{36})/a_{14}\}, \\
&& \{a_{12} = 0,  a_{13} = 0, a_{14} = 0, a_{26} = -(a_{15} a_{25})/a_{16}, a_{36} = -(a_{15} a_{35})/a_{16}, a_{46} = -(a_{15} a_{45})/a_{16}\},\\
&& \{a_{14} = 0, a_{15} = 0, a_{16} = 0,  a_{35} = -(a_{12} a_{25})/a_{13}, a_{36} = -(a_{12} a_{26})/a_{13}\},\\
&& \{a_{12} = 0, a_{13} = 0, a_{14} = 0, a_{15} = 0, a_{16} = 0\},\\
&& \{a_{13} = 0, a_{14} = 0, a_{15} = 0, a_{16} = 0, a_{25} = 0, a_{26} = 0\},\\
&& \{a_{12} = 0, a_{13} = 0, a_{14} = 0, a_{16} = 0, a_{25} = 0, a_{35} = 0, a_{45} = 0\}.
\end{eqnarray*}

From these solutions we obtain six classes of equigeodesic vectors, as shown below:  
\begin{center}
\begin{tabular}{l}
\thickline
\mbox{Equigeodesic vectors on} $\fr{m}_{12}\oplus\fr{m}_{23}$ \\
\hline
$ a_{12}\xi_{12}+a_{13}\xi_{13}+a_{14}\xi_{14}+a_{25}\xi_{25}+a_{26}\xi_{26}+a_{35}\xi_{35}+a_{36}\xi_{36} - \frac{a_{12}a_{25}+a_{13}a_{35}}{a_{14}}\xi_{45}- \frac{a_{12}a_{26}+a_{13}a_{36}}{a_{14}}\xi_{46}$, $a_{14}\neq 0$\\
$a_{12}\xi_{12} + a_{13}\xi_{13} + a_{25}\xi_{25} + a_{26}\xi_{26} + a_{45}\xi_{45} + a_{46}\xi_{46} - \frac{a_{12}a_{25}}{a_{13}}\xi_{35} - \frac{a_{12}a_{26}}{a_{13}}\xi_{36}$, $a_{13}\neq 0$\\
$a_{12}\xi_{12}+a_{35}\xi_{35}+a_{36}\xi_{36}+a_{45}\xi_{45}+a_{46}\xi_{46}$\\
\hline
\mbox{Equigeodesic vectors on} $\fr{m}_{13}\oplus\fr{m}_{23}$ \\
\hline
$a_{15}\xi_{15}+a_{16}\xi_{16}+a_{25}\xi_{25} - \frac{a_{15}a_{25}}{a_{16}}\xi_{26} + a_{35}\xi_{35} -\frac{a_{15}a_{35}}{a_{16}}\xi_{36} + a_{45}\xi_{45} - \frac{a_{15}a_{45}}{a_{16}}\xi_{46}$, $a_{16}\neq 0$\\
$a_{15}\xi_{15}+a_{26}\xi_{26}+a_{36}\xi_{36}+a_{46}\xi_{46}  $\\
\hline
\mbox{Equigeodesic vectors on} $\fr{m}_{23}$ \\
\hline
$a_{25}\xi_{25}+a_{26}\xi_{26}+ a_{35}\xi_{35}+a_{36}\xi_{36}+a_{45}\xi_{45}+a_{46}\xi_{46} $ (trivial)\\
\thickline
\end{tabular}
\end{center}

\begin{prop}
The generalized Wallach space $\SO(6)/\SO(3)\cdot\SO(2)$ admits six classes of equigeodesics vectors.  One of them is trivial, that is $X\in\fr{m}_{23}$.
\end{prop} 

For the special case where $n_1=n_2=1$ and $n_{3} = n-2$, that is the Stiefel manifold $V_{2}\bb{R}^{n} = \SO(n)/\SO(n-2)$, system (\ref{systemW}) becomes
\begin{eqnarray}\label{SolutionsV2Rn}
\left \{ \begin{array}{l}
 a_{12}a_{13} = a_{12}a_{14} = \cdots = a_{12}a_{1n} = 0\\
 a_{12}a_{23} = a_{12}a_{24} = \cdots = a_{12}a_{2n} = 0\\
 a_{13}a_{23} + a_{14}a_{24} + a_{15}a_{25} + \cdots + a_{1n-1}a_{2n-1} + a_{1n}a_{2n} = 0.
\end{array}\right.
\end{eqnarray} 

\noindent
$\blacktriangleright$ For $V_{2}\bb{R}^4$ the solutions of (\ref{SolutionsV2Rn}) are the following:
\begin{eqnarray*}
&&\{a_{12} = 0, a_{24} = -{a_{13} a_{23}}/{a_{14}}\}, \ \{a_{12} = 0, a_{13} = 0, a_{14} = 0\},\\
&& \{a_{12} = 0, a_{14} = 0, a_{23} = 0 \},  \{a_{13} = 0, a_{14} = 0, a_{23} = 0, a_{24} = 0\}.
\end{eqnarray*}
Therefore, the equigeodesic vectors  are
\begin{eqnarray*}
\left [ \begin{array}{ll}
 a_{13}\xi_{13} + a_{14}\xi_{14} + a_{23}\xi_{23} -\frac{a_{13}a_{23}}{a_{14}}\xi_{24} \ (a_{14}\neq 0), \  a_{23}\xi_{23} + a_{24}\xi_{24}, \ a_{13}\xi_{13} + a_{24}\xi_{24}, \ a_{12}\xi_{12}
\end{array} \right].
\end{eqnarray*}

\noindent
$\blacktriangleright$ For $V_{2}\bb{R}^5$ the solutions of (\ref{SolutionsV2Rn}) are
\begin{eqnarray*}
&&\{a_{12} = 0, a_{25} = -(a_{13} a_{23} + a_{14} a_{24})/a_{15} \}, \{a_{12} = 0, a_{15} = 0,  a_{24} = -a_{13} a_{23}/a_{14} \},\\
&&\{a_{12} = 0, a_{13} = 0, a_{14} = 0, a_{15} = 0\}, \{a_{12} = 0, a_{14} = 0, a_{15} = 0, a_{23} = 0\},\\
&& \{a_{13} = 0, a_{14} = 0, a_{15} = 0, a_{23} = 0, a_{24} = 0, a_{25} = 0\}.
\end{eqnarray*}
For this case we take the following classes of equigeodesic vectors:
\begin{eqnarray*}
\left [ \begin{array}{ll}
 a_{13}\xi_{13} + a_{14}\xi_{14} + a_{23}\xi_{23} -\frac{a_{13}a_{23}}{a_{14}}\xi_{24} + a_{25}\xi_{25}\ a_{14}\neq 0,  \  a_{23}\xi_{23} + a_{24}\xi_{24} + a_{25}\xi_{25}, \ a_{12}\xi_{12},\\
  a_{13}\xi_{13} + a_{24}\xi_{24}+ a_{25}\xi_{25}, \ 
 a_{13}\xi_{13} + a_{14}\xi_{14} + a_{15}\xi_{15} + a_{23}\xi_{23} +a_{24}\xi_{24} -\frac{(a_{13}a_{23}+a_{14}e_{24})}{a_{15}}\xi_{25},\ a_{15}\neq 0 
\end{array} \right ].
\end{eqnarray*}

\noindent
$\blacktriangleright$ For $V_{2}\bb{R}^6$ similarly we obtain the following classes equigeodesic vectors:
\begin{eqnarray*}
\left [ \begin{array}{ll}
 a_{13}\xi_{13} + a_{14}\xi_{14} + a_{15}\xi_{15} + a_{23}\xi_{23} + a_{24}\xi_{24} -\frac{(a_{13}a_{23}+ a_{14}e_{14})}{a_{15}}\xi_{25} + a_{26}\xi_{26},\ a_{15}\neq 0 \\ 
a_{23}\xi_{23} + a_{24}\xi_{24} + a_{25}\xi_{25} + a_{26}\xi_{26},\ a_{13}\xi_{13} + a_{24}\xi_{24}+ a_{25}\xi_{25} + a_{26}\xi_{26}, \ a_{12}\xi_{12},\\
 a_{13}\xi_{13} + a_{14}\xi_{14} + a_{23}\xi_{23} -\frac{a_{13}a_{23}}{a_{14}}\xi_{24} + a_{25}\xi_{25} + a_{26}\xi_{26},\ a_{14}\neq 0 \\
 a_{13}\xi_{13} + a_{14}\xi_{14} + a_{15}\xi_{15} + a_{16}\xi_{16} + a_{23}\xi_{23} +a_{24}\xi_{24} + a_{25}\xi_{25} -\frac{(a_{13}a_{23}+a_{14}\xi_{24}+ a_{15}\xi_{25})}{a_{16}}\xi_{26},\ a_{16}\neq 0
\end{array} \right ].
\end{eqnarray*}

\noindent
$\blacktriangleright$ For $V_{2}\bb{R}^n$ with $n\geq 7$ the equigeodesic vectors are the following:
\begin{eqnarray*}
\left [ \begin{array}{ll}
 \Span\{\xi_{23}, \xi_{24},\ldots, \xi_{2n}\}, \ \  \Span\{\xi_{12}\},\ \ a_{13}\xi_{13} + a_{24}\xi_{24} + \cdots  + a_{2n}\xi_{2n},\\
 a_{13}\xi_{13} + \cdots + a_{1n}\xi_{1n} + a_{23}\xi_{23} +\cdots + a_{2n-1}\xi_{2n-1} - \frac{(a_{13}a_{23} +\cdots + a_{1n-1}a_{2n-1})}{a_{1n}}\xi_{2n},\ a_{1n}\neq 0\\
 a_{13}\xi_{13} + \cdots + a_{1n-1}\xi_{1n-1} + a_{23}\xi_{23} +\cdots + a_{2n-2}\xi_{2n-2} - \frac{(a_{13}a_{23} +\cdots + a_{1n-2}a_{2n-2})}{a_{1n-1}}\xi_{2n-1}\\
  + a_{2n}\xi_{2n}, \ a_{1n-1}\neq 0\\
 a_{13}\xi_{13} + \cdots + a_{1n-2}\xi_{1n-2} + a_{23}\xi_{23} +\cdots + a_{2n-3}\xi_{2n-3} - \frac{(a_{13}a_{23} +\cdots + a_{1n-3}a_{2n-3})}{a_{1n-2}}\xi_{2n-2} 
 + a_{2n-1}\xi_{2n-1} \\+ a_{2n}\xi_{2n}, \ a_{1n-1}\neq 0\\
 a_{13}\xi_{13} + \cdots + a_{1n-3}\xi_{1n-3} + a_{23}\xi_{23} +\cdots + a_{2n-4}\xi_{2n-4} - \frac{(a_{13}a_{23} +\cdots + a_{1n-4}a_{2n-4})}{a_{1n-3}}\xi_{2n-3} \\
 + a_{2n-2}\xi_{2n-2} + a_{2n-1}\xi_{2n-1} + a_{2n}\xi_{2n}, \ a_{1n-3}\neq 0\\
 \quad \quad \quad \quad \quad \quad \quad \quad \quad \quad \quad \quad \quad \quad \quad \quad \quad \quad \quad \quad \vdots \\
 a_{13}\xi_{13} + a_{14}\xi_{14} + a_{15}\xi_{15} -\frac{(a_{13}a_{23} + a_{14}a_{24})}{a_{15}}\xi_{25} +\cdots + a_{2n}\xi_{2n},\ a_{15}\neq 0 \\
 a_{13}\xi_{13} + a_{14}\xi_{14} + a_{23}\xi_{23} -\frac{a_{13}a_{23}}{a_{14}}\xi_{24} +\cdots + a_{2n}\xi_{2n}, \ a_{14}\neq 0  
\end{array} \right ].
\end{eqnarray*}

Therefore, we obtain the following proposition
\begin{prop}
The Stiefel manifolds $V_{2}\bb{R}^n$ admit $n$ classes of equigeodesics vectors.  Two of them are trivial, that is $X\in\fr{m}_{12}$ and $X\in\fr{m}_{23}$.
\end{prop}

Next, we describe the equigeodesic vectors on $V_{2}\bb{R}^n$ with respect to the Einstein metric $\Lambda= (1, \lambda, \lambda)$ (cf. Theorem \ref{theoremEinstein}).  We call these vectors {\it Einstein equigeodesics vectors}.  

\begin{prop}
The Einstein manifold $(V_{2}\bb{R}^n, \Lambda = (1,\lambda,\lambda))$ admits two classes of Einstein equigeodesic vectors.  One is trivial and the other one belongs to the subspace $\fr{m}_{13}\oplus\fr{m}_{23}$.
\end{prop} 
\begin{proof}
Let $X\in\fr{m} = \sum_{1\leq i<j \leq 3}\fr{m}_{ij}$ as in (\ref{vectorW}).  Then equation $[X, \Lambda X]_{\fr{m}} = 0$ for $\Lambda = (1, \lambda, \lambda)$ is equivalent to the system
\begin{eqnarray*}
\left \{ \begin{array}{l}
 a_{12}a_{13} = a_{12}a_{14} = \cdots = a_{12}a_{1n} = 0\\
 a_{12}a_{23} = a_{12}a_{24} = \cdots = a_{12}a_{2n} = 0.
\end{array}\right.
\end{eqnarray*} 
This system has two solutions: $\{a_{12} = 0\}$ and $\{a_{13} = 0, a_{14} = 0,\ldots, a_{1n} = 0, a_{23} = 0, a_{24} = 0,\ldots, a_{2n} = 0\}$.  For the first solution the equigeodesic vectors are of the form $X = a_{13}\xi_{13} + a_{14}\xi_{14} + \cdots + a_{1n}\xi_{1n} + a_{23}\xi_{23} + a_{24}\xi_{24} + \cdots + a_{2n}\xi_{2n}$ that is $X\in\fr{m}_{13}\oplus\fr{m}_{23}$.  For the second solution the equigeodesic vectors belong to $\Span\{\xi_{12}\}$, i.e. they are trivial.
\end{proof}

\subsubsection{The Wallach space $W^6 = \U(3)/\U(1)^{3}$}
According to Subsection \ref{subsectionWallach} the tangent space of $W^6$ has three $\Ad^{\U(3)}$-invariant isotropy summands, $\fr{m}_{ij}$, $1\leq i<j \leq 3$ of dimensions $\dim\fr{m}_{ij} = 2$.  We set $\fr{m}_{ij} = \Span\{e_{ij}, f_{ij}\}$, where $e_{ij}, f_{ij}$ are the vectors as in (\ref{bracetSU}).  Let $X_{\fr{m}_{ij}} =a_{ij}e_{ij} + b_{ij}f_{ij}$ where $a_{ij}, b_{ij}\in\bb{R}$, be the vector on each irreducible component $\fr{m}_{ij}$ of the tangent space of the $W^6$.  Then we obtain the following: 

\begin{prop}\label{propositionW6}
We consider the Wallach space $W^6 = \U(3)/\U(1)^3$ and $X\in\fr{m}$.  Let $X = X_{\fr{m}_{12}}+X_{\fr{m}_{13}}+X_{\fr{m}_{23}}$, where $X_{\fr{m}_{ij}}$, $1\leq i<j \leq 3$ defined as above.  Then the equation $[X, \Lambda_{{\tiny\mathsf{Wal}}} X]_{\fr{m}} = 0$, for any invariant metric $\Lambda_{{\tiny\mathsf{Wal}}}$, is equivalent to the following system of algebraic equations:
\begin{eqnarray}\label{systemW6}
\left \{ \begin{array}{lll}
b_{12}b_{13} + a_{12}a_{13} =0, && b_{12}a_{13} - a_{12}b_{13} =0,\\
a_{12}a_{23} - b_{12}b_{23} =0, && a_{12}b_{23} + b_{12}a_{23} =0,\\
a_{13}a_{23} + b_{13}b_{23} =0, && a_{13}b_{23} - b_{13}a_{23} =0.
\end{array} \right.
\end{eqnarray}
Therefore, $X$ is an equigeodesic vector if and only if, $X\in\fr{m}_{ij}$, $(1\leq i< j\leq 3)$, that is $X$ is trivial.
\end{prop}
\begin{proof}
Let $X = \sum_{1\leq i< j\leq 3}X_{\fr{m}_{ij}}\in\fr{m}$.  Then according to  Lemmas \ref{BrackW}, \ref{basebracketU} and \ref{Lemma3} the vector $X$ is equigeodesic if and only if the coefficients of $X_{\fr{m}_{ij}}$ satisfy the system $\{[X_{\fr{m}_{12}},\, X_{\fr{m}_{13}}]_{\fr{m}} = 0, \, [X_{\fr{m}_{12}},\, X_{\fr{m}_{23}}]_{\fr{m}} = 0, \, [X_{\fr{m}_{13}},\, X_{\fr{m}_{23}}]_{\fr{m}} = 0\}$, that is the coefficients of $X$ satisfy the system of equations (\ref{systemW6}).  For the solutions of (\ref{systemW6}) we argue as follows:  Fix $i, j$, $1\leq i < j \leq 3$.

\noindent
\underline{Case I.}
If $a_{ij}, b_{ij}\neq 0$ and $a_{k\ell} = b_{k\ell} = 0$ with $a_{k\ell}\neq a_{ij}$ and $b_{k\ell}\neq  b_{ij}$, then the equigeodesic vectors are $X = a_{ij}e_{ij}+b_{ij}f_{ij}$.

\noindent
\underline{Case II.} 
If $a_{ij} = 0$ and $b_{ij} \neq0$ ($1\leq i < j \leq 3$), then the equigeodesic vectors are $X = b_{ij}f_{ij}$.

\noindent
\underline{Case III.}
If $a_{ij}\neq 0$ and $b_{ij} = 0$ ($1\leq i < j \leq 3$), then the equigeodesic vectors are $X = a_{ij}e_{ij}$.
\end{proof}

\begin{remark}
\textnormal{System (\ref{systemW6}) agrees with system of Example 3.8 of \cite{CGN} if we replace their blocks $\al_{ij}$, $1\leq i < j\leq 3$ with our variables $a_{ij} + ib_{ij}$ $(1\leq i < j \leq 3)$.}
\end{remark}

\subsubsection{The Wallach space $W^{12} = \Sp(3)/\Sp(1)^3$}
In this case the tangent space has three $\Ad^{\Sp(3)}$-invariant isotropy summands, $\fr{m}_{ij}$, $1\leq i < j\leq 3$ of dimensions $\dim\fr{m}_{ij} = 4$.  We set $\fr{m}_{ij} = \Span\{e_{ij}, f_{ij}, g_{ij},$ $h_{ij}\}$, where $e_{ij}, f_{ij},$ $g_{ij},$ $h_{ij}$ are the vectors as in (\ref{VectorsSP}).  Let $X_{\fr{m}_{ij}} =a_{ij}e_{ij} + b_{ij}f_{ij} + c_{ij}g_{ij} + q_{ij}h_{ij}$, where $a_{ij}, b_{ij}, c_{ij}, q_{ij}\in\bb{R}$, be the vector on each irreducible component $\fr{m}_{ij}$ of the tangent space of the $W^{12}$.  Then we have the following:

\begin{prop}\label{propositionW12}
We consider the Wallach space $\Sp(3)/\Sp(1)^3$ and $X\in\fr{m}$.  Let $X = X_{\fr{m}_{12}}+X_{\fr{m}_{13}}+X_{\fr{m}_{23}}$, where $X_{\fr{m}_{ij}}$, $1\leq i<j \leq 3$ defined as above.  Then the equation $[X, \Lambda_{{\tiny\mathsf{Wal}}} X]_{\fr{m}} = 0$ for any invariant metric $\Lambda_{{\tiny\mathsf{Wal}}}$ is equivalent to the following system of algebraic equations:
\begin{eqnarray}\label{systemW12}
\left \{ \begin{array}{lll}
a_{12}a_{13} + b_{12}b_{13} + c_{12}c_{13} + q_{12}q_{13} =0, && a_{12}b_{13} - b_{12}a_{13} + c_{12}q_{13} - q_{12}c_{13}=0,\\
a_{12}c_{13} - b_{12}q_{13} - c_{12}a_{13} + q_{12}b_{13} =0, && a_{12}q_{13} + b_{12}c_{13} - c_{12}b_{13} - q_{12}a_{13}=0,\\
a_{12}a_{23} - b_{12}b_{23} - c_{12}c_{23} - q_{12}q_{23} =0, && a_{12}b_{23} + b_{12}a_{23} - c_{12}q_{23} + q_{12}c_{23} =0, \\
a_{12}c_{23} + b_{12}q_{23} + c_{12}a_{23} - q_{12}b_{23} =0, && a_{12}q_{23} - b_{12}c_{23} + c_{12}b_{23} + q_{12}a_{23} =0,\\
a_{13}a_{23} + b_{13}b_{23} + c_{13}c_{23} + q_{13}q_{23} =0, && a_{13}b_{23} - b_{13}a_{23} - c_{13}q_{23} + q_{13}c_{23}= 0,\\
a_{13}c_{23} + b_{13}q_{23} - c_{13}a_{23} - q_{13}b_{23} =0, && a_{13}q_{23} - b_{13}c_{23} + c_{13}b_{23} - q_{13}a_{23} =0.
\end{array} \right.
\end{eqnarray}
Therefore, $X$ is an equigeodesic vector if and only if one of the following holds:
\begin{eqnarray*}
&& (1)\, X\in\fr{m}_{12}, \ (2)\, X\in\fr{m}_{13}, \ (3)\, X\in\fr{m}_{23},\ 
(4)\, X = a_{12}e_{12} + b_{12}f_{12} + c_{12}g_{12} - \frac{b_{12}c_{12}}{a_{12}}h_{12}\in\fr{m}_{12},\, a_{12}\neq 0\\
&& (5)\ X = a_{12}e_{12} + b_{12}f_{12} + c_{12}g_{12} + b_{23}f_{23} + c_{23}g_{23} + \frac{b_{12}c_{23}-b_{23}c_{12}}{a_{12}}h_{23}  \in\fr{m}_{12}\oplus\fr{m}_{23},\, a_{12}\neq 0.
\end{eqnarray*}
\end{prop}
\begin{proof}
Let $X_{\fr{m}_{ij}} =a_{ij}e_{ij} + b_{ij}f_{ij} + c_{ij}g_{ij} + q_{ij}h_{ij}\in\fr{m}_{ij}$.  Then according to the Lemmas \ref{BrackW}, \ref{LemmaLieSP} and \ref{Lemma3} the vector $X = \sum_{1\leq i< j} X_{ij}$ is equigeodesic if and only if the coefficients of $X_{\fr{m}_{ij}}$ satisfy the system $\{[X_{\fr{m}_{12}},\, X_{\fr{m}_{13}}]_{\fr{m}} = 0, \, [X_{\fr{m}_{12}},\, X_{\fr{m}_{23}}]_{\fr{m}} = 0, \, [X_{\fr{m}_{13}},\, X_{\fr{m}_{23}}]_{\fr{m}} = 0\}$, that is the coefficients of $X$ satisfy the system of equations (\ref{systemW12}). 
The solutions are the following:
\begin{eqnarray*}
&&\{a_{12} = 0, b_{12} = 0, c_{12} = 0, q_{12} = 0, a_{23} = 0, b_{23} = 0, c_{23} = 0, q_{23} = 0\},\\
&&\{a_{12} = 0, b_{12} = 0, c_{12} = 0, q_{12} = 0, a_{13} = 0, b_{13} = 0, c_{13} = 0, q_{13} = 0\},\\
&&\{a_{13} = 0, b_{13} = 0, c_{13} = 0, q_{13} = 0, a_{23} = 0, b_{23} = 0, c_{23} = 0, q_{23} = 0\},\\
&&\{a_{12} = 0, b_{12} = 0, c_{12} = 0, q_{12} = 0, a_{13} = 0, a_{23} = 0, b_{23} = 0, c_{23} = 0, q_{23} = 0\},\\
&&\{c_{12} = 0, a_{13} = 0, b_{13} = 0, c_{13} = 0, q_{13} = 0, a_{23} = 0, b_{23} = 0, c_{23} = 0, q_{23} = 0\},\\
&&\{a_{12} = 0, b_{12} = 0, a_{13} = 0, b_{13} = 0, c_{13} = 0, q_{13} = 0, a_{23} = 0, b_{23} = 0, c_{23} = 0, q_{23} = 0\},\\
&&\{a_{12} = 0, c_{12} = 0, a_{13} = 0, b_{13} = 0, c_{13} = 0, q_{13} = 0, a_{23} = 0, b_{23} = 0, c_{23} = 0, q_{23} = 0\},\\
&&\{a_{13} = 0, b_{13} = 0, c_{13} = 0, q_{13} = 0, a_{23} = 0, b_{23} = 0, c_{23} = 0, q_{23} = 0, q_{12} = -\frac{b_{12}c_{12}}{a_{12}}\},\\
&&\{a_{13} = 0, b_{13} = 0, c_{13} = 0, q_{13} = 0, a_{23} = 0, q_{23} = \frac{b_{12}c_{23}-b_{23}c_{12}}{a_{12}}\}.
\end{eqnarray*}
Therefore, we determine all vectors that appear in Proposition \ref{propositionW12}.
\end{proof}

\subsection{Algebraic equigeodesic vectors on the Stiefel manifolds $V_{1+n_2}\bb{R}^n$} 

When the Stiefel manifold $G/H = V_{1+n_2}\bb{R}^{n}$ is a total space over the generalized Wallach space $P = \SO(n)/\SO(n_2)\cdot\SO(n_3)$ the tangent space $T_{o}(G/H)$ is written as a direct sum of the tangent space of the base space $T_{o}P = \fr{m}_{12}\oplus\fr{m}_{13}\oplus\fr{m}_{23}$ and the tangent space of the fiber $F = \SO(n_2)$, that is $T_oF = \fr{so}(n_2)$.   A basis for $\fr{m}_{ij}$ $(i\leq i< j\leq 3)$ is given in (\ref{baseWallach}) and for $\fr{so}(n_2)$ it is $\Span\{\xi_{ab} : 2\leq a\leq 1+n_2, 3\leq b\leq 1+n_2\}$.  For each $X\in\fr{so}(n_2)$, write $X_{\fr{so}(n_2)} = \sum_{i=2}^{1+n_2}\sum_{j = 3}^{1+n_2} a_{ij}\xi_{ij}$.  Then by the same way as in Proposition \ref{PropWallach} we have:  

\begin{prop}
We consider the Stiefel manifold $V_{1+n_2}\bb{R}^n$ and for $X\in\fr{m}$ we write $X = X_{\fr{so}(n_2)} + \sum_{1\leq i<j<3}X_{\fr{m}_{ij}}$, where $X_{\fr{m}_{ij}}$ are as in {\rm (\ref{vectorW})} and $X_{\fr{so}(n_2)}$ as above.  Then the vector $X$ is an euigeodesic vector if and only if  
{\small \begin{eqnarray}\label{systemSt}
{\mathsf{(A)}}\left \{ \begin{array}{l} 
\sum_{i\in I_1, j\in J_2}\sum_{k\in I_1, l\in J_3} a_{ij}a_{kl} B([\xi_{ij}, \xi_{ab}]_{\fr{m}}, \xi_{kl})=0\\ 
 \quad \quad \quad \quad \quad \quad \quad \quad  \    \vdots\\ 
\sum_{i\in I_1, j\in J_2}\sum_{k\in I_1, l\in J_3} a_{ij}a_{kl} B([\xi_{ij}, \xi_{d_3}]_{\fr{m}}, \xi_{kl})=0\nonumber\\ 
\end{array} \right. &
{\mathsf{(D)}}\left \{ \begin{array}{l} 
\sum_{i\in I_1, j\in J_2}\sum_{k\in K, l\in L} a_{ij}a_{kl} B([\xi_{ij}, \xi_{ab}]_{\fr{m}}, \xi_{kl})=0\\ 
 \quad \quad \quad \quad \quad \quad \quad \quad  \    \vdots\\ 
\sum_{i\in I_1, j\in J_2}\sum_{k\in K, l\in L} a_{ij}a_{kl} B([\xi_{ij}, \xi_{d_1}]_{\fr{m}}, \xi_{kl})=0\nonumber\\ 
\end{array} \right.\\
{\mathsf{(B)}}\left \{ \begin{array}{l} 
\sum_{i\in I_1, j\in J_2}\sum_{k\in I_2, l\in J_3} a_{ij}a_{kl} B([\xi_{ij}, \xi_{ab}]_{\fr{m}}, \xi_{kl})=0\\ 
 \quad \quad \quad \quad \quad \quad \quad \quad  \    \vdots \\ 
\sum_{i\in I_1, j\in J_2}\sum_{k\in I_2, l\in J_3} a_{ij}a_{kl} B([\xi_{ij}, \xi_{d_2}]_{\fr{m}}, \xi_{kl})=0\\ 
\end{array} \right. &
{\mathsf{(E)}}\left \{ \begin{array}{l} 
\sum_{i\in I_2, j\in J_3}\sum_{k\in K, l\in L} a_{ij}a_{kl} B([\xi_{ij}, \xi_{ab}]_{\fr{m}}, \xi_{kl})=0\\ 
 \quad \quad \quad \quad \quad \quad \quad \quad  \    \vdots\\ 
\sum_{i\in I_2, j\in J_3}\sum_{k\in K, l\in L} a_{ij}a_{kl} B([\xi_{ij}, \xi_{d_3}]_{\fr{m}}, \xi_{kl})=0\nonumber\\ 
\end{array} \right.\\
{\mathsf{(C)}}\left \{ \begin{array}{l} 
\sum_{i\in I_1, j\in J_3}\sum_{k\in I_2, l\in J_3} a_{ij}a_{kl} B([\xi_{ij}, \xi_{ab}]_{\fr{m}}, \xi_{kl})=0\\ 
 \quad \quad \quad \quad \quad \quad \quad \quad \      \vdots\\ 
\sum_{i\in I_1, j\in J_3}\sum_{k\in I_2, l\in J_3} a_{ij}a_{kl} B([\xi_{ij}, \xi_{d_1}]_{\fr{m}}, \xi_{kl})=0.\nonumber 
\end{array} \right. &
\end{eqnarray} }
where $K = \{k : 2 \leq k \leq n_2\}$ and $L = \{l : 3 \leq l \leq 1+n_2\}$.






\end{prop}

From the above Proposition it is easy to see that for $n_2=3$, $n_3 = 2$ and for the Jensen Einstein metrics (c.f. Remark \ref{remarkJensen}) we obtain
\begin{prop}
The Stiefel manifold $(V_{4}\bb{R}^6, \Lambda_{\tiny\mathsf{Jen}} = (\lambda_2, \lambda_2, 1, 1))$ admits $42$ classes of Jensen's Einstein equigeodesics vectors. 
\end{prop}
\begin{proof}
Let $X_{\fr{m}_{ij}} = \sum_{k\in I, l\in J} a_{kl} \xi_{kl}\in\fr{m}_{ij}$ and $X_{\fr{so}(3)} = \sum_{2\leq i < j\leq4}a_{ij}\xi_{ij}$.  Then the vector $X = X_{\fr{so}(3)} + X_{\fr{m}_{12}} + X_{\fr{m}_{13}} + X_{\fr{m}_{23}}$ is an equigeodesic vector with respect to  Jensen's Einstein metric $\Lambda_{\tiny\mathsf{Jen}}$,  if and only if it satisfies the system $\{[X_{\fr{m}_{23}}, X_{\fr{so}(3)}]_{\fr{m}} = 0, \ [X_{\fr{m}_{13}}, X_{\fr{m}_{12}}]_{\fr{m}} = 0, \ [X_{\fr{m}_{23}}, X_{\fr{m}_{12}}]_{\fr{m}} = 0\}$.  According to Lemma \ref{bracSO} we have $[\fr{so}(3), \fr{m}_{23}]\subset\fr{m}_{23}$, therefore the previous system is equivalent to 
\begin{eqnarray*}
\left \{ \begin{array}{lll}
a_{15}a_{12} - a_{45}a_{24} - a_{35}a_{23} = 0, && a_{16}a_{12} - a_{46}a_{24} - a_{36}a_{23} = 0,\\
a_{25}a_{23} + a_{45}a_{34} + a_{15}a_{13} = 0, && a_{26}a_{23} - a_{46}a_{34} + a_{16}a_{13} = 0,\\
a_{25}a_{24} + a_{35}a_{34} + a_{15}a_{14} = 0, && a_{26}a_{24} + a_{36}a_{34} + a_{16}a_{14} = 0,\\
a_{25}a_{12} + a_{35}a_{13} + a_{45}a_{14} = 0, && a_{26}a_{12} + a_{36}a_{13} + a_{46}a_{14} = 0.
\end{array} \right.
\end{eqnarray*}
The solutions are given in {\bf Table 1.}
\end{proof}

From {\bf Table 1.} it is easy to see that some of the 42 classes of equigeodesic vectors are the following:
\begin{eqnarray*}
\begin{array}{ll} 
{\bf 1.} & \frac{a_{14}a_{23}}{a_{24}}\xi_{13} + a_{14}\xi_{14} + a_{15}\xi_{15} + a_{16}\xi_{16} + a_{23}\xi_{23} + a_{24}\xi_{24} -\frac{a_{14}a_{14}}{a_{24}}\xi_{25} - \frac{a_{14}a_{16}}{a_{24}}\xi_{26} + a_{34}\xi_{34} + a_{35}\xi_{35} \\
&+ a_{36}\xi_{36} + \frac{a_{12}a_{15} - a_{23}a_{35}}{a_{24}}\xi_{45} + \frac{a_{12}a_{16} - a_{23}a_{36}}{a_{24}}\xi_{46}, a_{24}\neq 0\\
{\bf 2.}&  a_{12}\xi_{12} + a_{13}\xi_{13} + a_{14}\xi_{14} + a_{25}\xi_{25} + a_{26}\xi_{26} + a_{35}\xi_{35} + a_{36}\xi_{36} - \frac{a_{12}a_{25}+a_{13}a_{35}}{a_{14}}\xi_{45}  \\
& - \frac{a_{12}a_{26} + a_{13}a_{36}}{a_{14}}\xi_{46}, a_{14}\neq 0\\
{\bf 3.}&  a_{23}\xi_{23} + a_{24}\xi_{24} + a_{25}\xi_{25} + a_{26}\xi_{26} + a_{34}\xi_{34} - \frac{a_{24}a_{25}}{a_{34}}\xi_{35} - \frac{a_{24}a_{26}}{a_{34}}\xi_{36} + \frac{a_{23}a_{25}}{a_{34}}\xi_{45} + \frac{a_{23}a_{26}}{a_{34}}\xi_{46}, a_{34}\neq 0\\
{\bf 4.}&  a_{13}\xi_{13} + a_{14}\xi_{14} + a_{15}\xi_{15} + a_{16}\xi_{16} +a_{34}\xi_{34} - \frac{a_{14}a_{15}}{a_{34}}\xi_{35} - \frac{a_{14}a_{16}}{a_{34}}\xi_{36} + \frac{a_{13}a_{15}}{a_{34}}\xi_{45} + \frac{a_{13}a_{16}}{a_{34}}\xi_{46}, a_{34}\neq 0\\
{\bf 5.}& a_{12}\xi_{12} + a_{14}\xi_{14} + a_{15}\xi_{15} + a_{16}\xi_{16} + a_{24}\xi_{24} - \frac{a_{14}a_{15}}{a_{24}}\xi_{25} - \frac{a_{14}a_{16}}{a_{24}}\xi_{26} + \frac{a_{12}a_{15}}{a_{24}}\xi_{45} + \frac{a_{12}a_{16}}{a_{24}}\xi_{46}, a_{24}\neq 0\\
{\bf 6.} & a_{12}\xi_{12} + a_{16}\xi_{16} + a_{23}\xi_{23} + a_{24}\xi_{24} + a_{35}\xi_{35} + a_{36}\xi_{36} - \frac{a_{23}a_{35}}{a_{24}}\xi_{45} + \frac{a_{12}a_{16} - a_{23}a_{36}}{a_{24}}\xi_{46}, a_{24}\neq 0\\
{\bf 7.}&  a_{23}\xi_{23} + a_{24}\xi_{24} + a_{35}\xi_{35} + a_{36}\xi_{36} -\frac{a_{23}a_{35}}{a_{24}}\xi_{45} - \frac{a_{23}a_{36}}{a_{24}}\xi_{46}, a_{24}\neq 0\\
{\bf 8.}&  a_{12}\xi_{12} + a_{15}\xi_{15} + a_{16}\xi_{16} + a_{23}\xi_{23} + \frac{a_{12}a_{15}}{a_{23}}\xi_{35} + \frac{a_{12}a_{16}}{a_{23}}\xi_{36}, a_{23}\neq 0.
\end{array}
\end{eqnarray*}

\begin{landscape}
\begin{center}
{{\bf Table 1.}  Jensen's Einstein equigeodesic vectors for the Stiefel manifold} $V_{4}\bb{R}^6$ with tangent space $\fr{so}(3)\oplus\fr{m}_{12}\oplus\fr{m}_{13}\oplus\fr{m}_{23}$
\end{center}
\begin{center}
\begin{tabular}{ll}
\thickline
$\{a_{12} =\frac{-a_{14} a_{23} + a_{13} a_{24}}{a_{34}}, a_{35} = \frac{-a_{14} a_{15} - a_{24} a_{25}}{a_{34}}, a_{36} = \frac{-a_{14} a_{16} - a_{24} a_{26}}{a_{34}}, a_{45} = \frac{ a_{13} a_{15} + a_{23} a_{25}}{a_{34}}, a_{46} = \frac{a_{13} a_{16} + a_{23} a_{26}}{a_{34}}\} $
\\
$\{a_{12} = 0, a_{13} = \frac{a_{14} a_{23}}{a_{24}}, a_{35} = \frac{-a_{14} a_{15} - a_{24} a_{25}}{a_{34}}, a_{36} = \frac{-a_{14} a_{16} -  a_{24} a_{26}}{a_{34}}, a_{45} = \frac{a_{23} (a_{14} a_{15} + a_{24} a_{25})}{a_{24} a_{34}}, a_{46} = \frac{ a_{23} (a_{14} a_{16} + a_{24} a_{26})}{a_{24} a_{34}}\} $
\\
$\{a_{12} = \frac{a_{13} a_{24}}{a_{34}}, a_{14} = 0, a_{35} = \frac{-a_{24} a_{25}}{a_{34}}, a_{36} = \frac{-a_{24} a_{26}}{a_{34}}, a_{45} = \frac{a_{13} a_{15} + a_{23} a_{25}}{a_{34}}, a_{46} = \frac{a_{13} a_{16} + a_{23} a_{26}}{a_{34}}\}
$
\\
$\{a_{13} = \frac{a_{14} a_{23}}{a_{24}}, a_{34} = 0, a_{25} = \frac{-a_{14} a_{15}}{a_{24}}, a_{26} = \frac{-a_{14} a_{16}}{a_{24}}, a_{45} = \frac{a_{12} a_{15} - a_{23} a_{35}}{a_{24}}, a_{46} = \frac{
 a_{12} a_{16} - a_{23} a_{36}}{a_{24}}\}
$
\\
$\{a_{12} = 0, a_{13} = \frac{a_{14} a_{23}}{a_{24}}, a_{25} = \frac{-a_{14} a_{15}}{a_{24}}, a_{35} = 0, a_{36} = \frac{-a_{14} a_{16} - a_{24} a_{26}}{a_{34}}, a_{45} = 0, a_{46} = \frac{a_{23} (a_{14} a_{16} + a_{24} a_{26})}{a_{24} a_{34}}\}
$
\\
$\{a_{12} = 0, a_{13} = 0, a_{14} = 0, a_{35} = \frac{-a_{24} a_{25}}{a_{34}}, a_{36} = \frac{-a_{24} a_{26}}{a_{34}}, a_{45} = \frac{a_{23} a_{25}}{a_{34}}, a_{46} = \frac{a_{23} a_{26}}{a_{34}}\}
$
\\
$\{a_{12} = 0, a_{13} = 0, a_{23} = 0, a_{35} = \frac{(-a_{14} a_{15} - a_{24} a_{25})}{a_{34}}, a_{36} = \frac{(-a_{14} a_{16} - a_{24} a_{26})}{a_{34}}, a_{45} = 0, a_{46} = 0\}
$
\\
$\{a_{12} =  a_{14} =  a_{24} =  a_{35} =  a_{36} = 0, a_{45} =\frac{a_{13} a_{15} + a_{23} a_{25}}{a_{34}}, a_{46} = \frac{a_{13} a_{16} + a_{23} a_{26}}{a_{34}}\}$
\\
$\{a_{12} =  a_{23} = a_{24} = 0, a_{35} =\frac{-a_{14} a_{15}}{a_{34}}, a_{36} = \frac{-a_{14} a_{16}}{a_{34}}, a_{45} = \frac{a_{13} a_{15}}{a_{34}}, a_{46} = \frac{a_{13} a_{16}}{a_{34}}\}
$
\\
$\{a_{12} = 0, a_{13} = \frac{a_{14} a_{23}}{a_{24}}, a_{34} = 0, a_{25} =\frac{-a_{14} a_{15}}{a_{24}}, a_{26} = \frac{-a_{14} a_{16}}{a_{24}}, a_{45} = \frac{-a_{23} a_{35}}{a_{24}},  a_{46} = \frac{-a_{23} a_{36}}{a_{24}} \}
$
\\
$\{a_{13} = \frac{a_{14} a_{23}}{a_{24}}, a_{15} = 0, a_{34} = a{25} = 0, a_{26} = \frac{- a_{14} a_{16}}{a_{24}}, a_{45} = \frac{-a_{23} a_{35}}{a_{24}}, a_{46} = \frac{a_{12} a_{16} - a_{23} a_{36}}{a_{24}}\} $
\\
$\{a_{13} = a_{14} = a_{34} = a_{25} = a_{26} = 0, a_{45} = \frac{a_{12} a_{15} - a_{23} a_{35}}{a_{24}}, a_{46} = \frac{a_{12} a_{16} - a_{23} a_{36}}{a_{24}} \}
$
$\{a_{13} = a_{23} = a_{34} = 0, a_{25} = \frac{-a_{14} a_{15}}{a_{24}}, a_{26} = \frac{-a_{14} a_{16}}{a_{24}}, a_{45} = \frac{a_{12} a_{15}}{a_{24}}, a_{46} = \frac{a_{12} a_{16}}{a_{24}} \}
$
\\
$\{a_{14} = a_{24} = a_{34} = 0, a_{25} = \frac{-a_{13} a_{15}}{a_{23}}, a_{26} = \frac{-a_{13} a_{16}}{a_{23}}, a_{35} = \frac{a_{12} a_{15}}{a_{23}}, a_{36} = \frac{a_{12} a_{16}}{a_{23}} \}
$
\\
$\{a_{15} = a_{16} = a_{23} = a_{24} = a_{34} = 0, a_{45} = \frac{-a_{12} a_{25} - a_{13} a_{35}}{a_{14}}, a_{46} = \frac{-a_{12} a_{26} - a_{13} a_{36}}{a_{14}} \} 
$
$\{a_{12} = a_{13} = a_{14} = a_{25} = a_{35} = 0, a_{36} = \frac{-a_{24} a_{26}}{a_{34}}, a_{45} = 0, a_{46} = \frac{a_{23} a_{26}}{a_{34}} \}
$
\\
$\{a_{12} = a_{13} = a_{23} = 0, a_{25} = \frac{-a_{14} a_{15}}{a_{24}}, a_{35} = 0, a_{36} = \frac{-a_{14} a_{16} - a_{24} a_{26}}{a_{34}}, a_{45} = a_{46} = 0 \} 
$
\\
$\{a_{12} = a_{14} = a_{24} = 0, a_{25} = \frac{-a_{13} a_{15}}{a_{23}}, a_{35} = a_{36} = a_{45} = 0, a_{46} = \frac{a_{13} a_{16} + a_{23} a_{26}}{a_{34}} \}
$
$\{a_{12} = a_{15} = a_{23} = a_{24} = a_{35} = 0, a_{36} = \frac{-a_{14} a_{16}}{a_{34}}, a_{45} = 0, a_{46} = \frac{a_{13} a_{16}}{a_{34}} \}
$
\\
$\{a_{12} = 0, a_{13} = \frac{a_{14} a_{23}}{a_{24}}, a_{15} = a_{34} =  a_{25} = 0, a_{26} = \frac{-a_{14} a_{16}}{a_{24}}, a_{45} = \frac{-a_{23} a_{35}}{a_{24}}, a_{46} = \frac{-a_{23} a_{36}}{a_{24}} \} 
$
\\
$\{a_{12} = a_{13} = a_{14} = a_{34} = a_{25} = a_{26} = 0, a_{45} = \frac{-a_{23} a_{35}}{a_{24}}, a_{46} = \frac{-a_{23} a_{36}}{a_{24}} \} 
$
$\{a_{13} = a_{14} = a_{15} = a_{34} = a_{25} = a_{26} = 0, a_{45} = \frac{-a_{23} a_{35}}{a_{24}}, a_{46} = \frac{a_{12} a_{16} - a_{23} a_{36}}{a_{24}} \}
$
\\
$\{a_{13} = a_{14} = a_{23} = a_{34} = a_{25} = a_{26} = 0, a_{45} = \frac{a_{12} a_{15}}{a_{24}}, a_{46} = \frac{a_{12} a_{16}}{a_{24}} \}
$
$\{a_{13} = a_{15} = a_{23} = a_{34} = a_{25} = 0, a_{26} = \frac{-a_{14} a_{16}}{a_{24}}, a_{45} = 0, a_{46} = \frac{a_{12} a_{16}}{a_{24}} \}
$
\\
$\{a_{12} = a_{14} = a_{24} = a_{34} = 0, a_{25} = \frac{-a_{13} a_{15}}{a_{23}}, a_{26} = \frac{-a_{13} a_{16}}{a_{23}}, a_{35}  = a_{36} = 0 \}
$
$\{a_{13} = a_{14} = a_{24} = a_{34} = a_{25} = a_{26} = 0, a_{35} = \frac{a_{12} a_{15}}{a_{23}}, a_{36} = \frac{a_{12} a_{16}}{a_{23}} \}
$
\\
$\{a_{14} = a_{15} = a_{24} = a_{34} = a_{25} = 0, a_{26} = \frac{-a_{13} a_{16}}{a_{23}}, a_{35} = 0, a_{36} = \frac{a_{12} a_{16}}{a_{23}} \}
$
$\{a_{14} = a_{15} = a_{16} = a_{23} = a_{24} = a_{34} = 0, a_{35} = \frac{-a_{12} a_{25}}{a_{13}}, a_{36} = \frac{-a_{12} a_{26}}{a_{13}} \} 
$
\\
$\{a_{13} = \frac{a_{14} a_{23}}{a_{24}}, a_{15} = a_{16} = a_{25} = a_{26} = a_{35} = a_{36} = a_{45} = a_{46} = 0 \}
$
$\{a_{12} = a_{14} = a_{15} = a_{23} = a_{24} = a_{35} = a_{36} = a_{45} = 0, a_{46} = \frac{a_{13} a_{16}}{a_{34}} \}
$
\\
$\{a_{12} = a_{13} = a_{14} = a_{15} = a_{34} = a_{25} = a_{26} = 0, a_{45} = \frac{-a_{23} a_{35}}{a_{24}}, a_{46} = \frac{-a_{23} a_{36}}{a_{24}} \}
$
$\{a_{13} = a_{14} = a_{15} = a_{23} = a_{34} = a_{25} = a_{26} = a_{45} = 0, a_{46} = \frac{a_{12} a_{16}}{a_{24}} \}
$
\\
$\{a_{12} = a_{14} = a_{15} = a_{24} = a_{34} = a_{25} = 0, a_{26} = \frac{-a_{13} a_{16}}{a_{23}}, a_{35} = a_{36} = 0 \}
$
$\{a_{13} = a_{14} = a_{15} = a_{24} = a_{34} = a_{25} = a_{26} = a_{35} = 0, a_{36} = \frac{a_{12} a_{16}}{a_{23}} \}
$
\\
$\{a_{12} = a_{14} = a_{15} = a_{16} = a_{23} = a_{24} = a_{34} = a_{35} = a_{36} = 0 \}
$
$\{a_{13} = a_{14} = a_{15} = a_{16} = a_{23} = a_{24} = a_{34} = a_{25} = a_{26} = 0 \}
$
\\
$\{a_{14} = a_{15} = a_{16} = a_{24} = a_{25} = a_{26} = a_{35} = a_{36} = a_{45} = a46 = 0 \}
$
$\{a_{15} = a_{16} = a_{23} = a_{24} = a_{25} = a_{26} = a_{35} = a_{36} = a_{45} = a_{46} = 0 \}\
$
\\
$\{a_{12} = a_{13} = a_{14} = a_{23} = a_{24} = a_{35} = a_{36} = a_{45} = a_{46} = 0 \} 
$
$\{a_{14} = a_{15} = a_{16} = a_{25} = a_{26} = a_{35} = a_{36} = a_{45} = a_{46} = 0 \}
$
\\
$\{a_{15} = a_{16} = a_{25} = a_{26} = a_{35} = a_{36} = a_{45} = a_{46} = 0 \}
$
$\{a_{12} = a_{13} = a_{14} = a_{23} = a_{24} = a_{34} = 0 \}
$
\\
\thickline
\end{tabular}
\end{center}
\end{landscape}

\subsection{Algebraic equigeodesic vectors on spheres} 

\subsubsection{The sphere $\bb{S}^{2n+1} = \U(n+1)/\U(n)$}
As we saw in Section \ref{sectionSphere1}, the tangent space of $\U(n+1)/\U(n)$ is written as a direct sum of two irreducible modules of dimensions 1 and $2n$.  In particular, we take $\fr{m}_1 = \Span\{F_{11} = 1/2f_{11}\}$ and $\fr{m}_2 = \Span\{e_{1j}, f_{1j} : 1 < j \leq n+1\}$.  
Then on $\bb{S}^{2n+1}$ the equigeodesic vectors are all trivial.  In particular, we deduce the following:
\begin{prop}
We consider the sphere $\U(n+1)/\U(n)$.  Let $X\in\fr{m} \cong T_{e\U(n+1)}\bb{S}^{2n+1}$ be a non zero vector.  Then $X$ is an equigeodesic vector if and only if $X\in\fr{m}_1$ or $X\in\fr{m}_2$. 
\end{prop} 
\begin{proof}
Let $X = X_{\fr{m}_1} + X_{\fr{m}_2}$, where $X_{\fr{m}_1} = \al_{11} F_{11}$, $\al_{11}\in\bb{R}$ and $X_{\fr{m}_2} = \sum_{j = 2}^{n+1} \al_{1j}e_{1j} + \beta_{1j}f_{ij}$, $\al_{1j}, \beta_{1j}\in \bb{R}$.  Then by  Lemma \ref{Lemma3}, $X$ is an equigeodesic vector if  and only if $[X_{\fr{m}_1}, \, X_{\fr{m}_2}]_{\fr{m}} = 0$.  Therefore we have that
\begin{eqnarray*}
&& 0 = [X_{\fr{m}_1}, \, X_{\fr{m}_2}]_{\fr{m}} = [\al_{11} F_{11} \, ,\, \sum_{j = 2}^{n+1} (\al_{1j}e_{1, j} + \beta_{1j}f_{ij})]_{\fr{m}} = \al_{11} ( \sum_{j = 2}^{n+1} (\al_{1j}[F_{11}, e_{1j}]_{\fr{m}} + \beta_{1j}[F_{11}, f_{1j}]_{\fr{m}})).
\end{eqnarray*}
Hence we conclude that,  if $\al_{11} = 0$ then $\al_{1j}, \beta_{1j}\in \bb{R}$, for $j = 2,3,\ldots, n+1$, that is $X\in\fr{m}_2$ and  if $\al_{11}\in\bb{R}$ then $\al_{12} = \beta_{12} = \al_{13} = \cdots = \al_{1n+1} = \beta_{1n+1} = 0$, that is $X\in\fr{m}_1$.
\end{proof}

\subsubsection{The sphere $\bb{S}^{4n+3} = \Sp(n+1)/\Sp(n)$}
Recall from \ref{sectionSphere2} that the tangent space of $\bb{S}^{4n+3}$ splits into two real irreducible $\Ad^{\Sp(n)}$-invariant summands, that is $\fr{m} = \fr{m}_1\oplus\fr{m}_2$.  We take $\fr{m}_1 = \Span\{u_{11}, k_{1n+2},$ $\mu_{1n+2}\}$ and $\fr{m}_2 = \Span\{e_{12}, f_{12},\ldots e_{1n+1}, f_{1n+1}, g_{12}, h_{12},\ldots, g_{1n+1}, h_{1n+1}\}$.  Let $X = X_{\fr{m}_1} + X_{\fr{m}_2}$, where $X_{\fr{m}_1} = a_{11}u_{11} + a_{1n+2}k_{1n+2} + b_{1n+2}\mu_{1n+2}$ with $a_{11}, a_{1n+2}, b_{1n+2}\in\bb{R}$ and $X_{\fr{m}_2} = \sum_{i=2}^{n+1}c_{1i}e_{1i} + d_{1i}f_{1i} + l_{1i}g_{1i} + m_{1i}h_{1i}$, with $c_{1i}, d_{1i}, l_{1i}, m_{1i}\in\bb{R}$.  Then, a vector $X\in\fr{m}$ it is equigeodesic if and only if $[X_{\fr{m}_1}, X_{\fr{m}_2}]_{\fr{m}} = 0$.  According to  Lemmas \ref{LemmaLieSP} and \ref{LemmaLieSP2} we obtain the following proposition:

\begin{prop}
We consider the sphere $\Sp(n+1)/\Sp(n)$ and $X\in\fr{m}$.  Let $X = X_{\fr{m}_1}+X_{\fr{m}_2}$ where $X_{\fr{m}_1}$, $X_{\fr{m}_2}$ as above.  Then the equation $[X, \Lambda X]_{\fr{m}} = 0$ for any invariant metric $\Lambda$ is equivalent to the following system of algebraic equations:
\begin{eqnarray}\label{sphereSP}
\left \{ \begin{array}{l} 
a_{11}d_{12}+a_{1n+2}l_{12}+b_{1n+2}m_{12} = 0\\
a_{11}c_{12}-a_{1n+2}m_{12} +b_{1n+2}l_{12} = 0\\
\quad \quad \ \quad \ \quad \ \quad \ \quad \vdots\\
a_{11}d_{n+1}+a_{1n+2}l_{1n+1}+b_{1n+2}m_{1n+1} = 0\\
a_{11}c_{1n+1}-a_{1n+2}m_{1n+1}+b_{1n+2}l_{1n+1} = 0\\
a_{1n+2}c_{12} - b_{1n+2}d_{12} + a_{11}m_{12} = 0\\
-a_{11}l_{12} + a_{1n+2}d_{12} + c_{12}b_{1n+2} = 0\\
\quad \quad \ \quad \ \quad \ \quad \ \quad \vdots\\
a_{11}m_{1n+1} + a_{1n+2}c_{1n+1} - b_{1n+2}d_{1n+1} = 0\\
-a_{11}l_{1n+1} + a_{1n+2}d_{1n+1} + b_{1n+2}c_{1n+1} = 0.
\end{array} \right. 
\end{eqnarray}
The solutions of this system are  
\begin{eqnarray*}
&& \mathcal{A} := a_{11} =  a_{1n+2} =  b_{1n+2} = 0 \\
&& \mathcal{B} := c_{12} = d_{12} =\cdots = c_{1n+1} = d_{1n+1} = l_{12} = m_{12} = \cdots = l_{1n+1} = m_{1n+1} = 0\\
&& \mathcal{C} := a_{11} = c_{12} = d_{12} =\cdots = c_{1n+1} = d_{1n+1} = l_{12} = m_{12} = \cdots = l_{1n+1} = m_{1n+1} = 0 
\end{eqnarray*}
and all the combinations $a_{11} = a_{1n+2} = \mathcal{B}$, $a_{11} = b_{1n+2} =  \mathcal{B}$, $a_{1n+2} =  \mathcal{B}$, $b_{1n+2} =  \mathcal{B}$, $a_{1n+2} = b_{1n+2} = \mathcal{B}$, $c_{1n+1} =  \mathcal{A}$, $d_{1n+1} =  \mathcal{A}$, $l_{1n+1} = \mathcal{A}$, $m_{1n+1} =  \mathcal{A}$, $c_{1n+1} = d_{1n+1} = \mathcal{A}$, $c_{1n+1} = d_{1n+1} = l_{1n+1} =  \mathcal{A}$, $a_{1n+2} = c_{1n+1} = d_{1n+1} =  \mathcal{B}$, etc.  Therefore, $X$ is an equigeodesic vector if and only if $X\in\fr{m}_1$ or $X\in\fr{m}_2$.
\end{prop}

\begin{example}\textnormal{
For $\bb{S}^{7} = \Sp(2)/\Sp(1)$ system (\ref{sphereSP}) takes the form
\begin{eqnarray*}
\left \{ \begin{array}{ll}
a_{11}d_{12} + a_{13}l_{12} + b_{13}m_{12} = 0 & a_{11}c_{12} - a_{13}m_{12} + b_{13}l_{12} = 0\\
a_{11}m_{12} + a_{13}c_{12} - b_{13}d_{12} = 0 & -a_{11}l_{12} + a_{13}d_{12} + b_{13}c_{12} = 0.
\end{array}\right.
\end{eqnarray*} 
Some of the solutions are
$
\{a_{11} = 0, a_{13} = 0, b_{13} = 0\},\, \{c_{12} = 0, d_{12} = 0, l_{12} = 0, m_{12} = 0\},\, \{a_{11} = 0,$ $c_{12} = 0, d_{12} = 0, l_{12} = 0, m_{12} = 0\}$, $\{a_{11} = 0, a_{13} = 0, c_{12} = 0, d_{12} = 0, l_{12} = 0, m_{12} = 0\}$, $\{b_{13} = 0, c_{12} = 0, d_{12} = 0, l_{12} = 0, m_{12} = 0\}$, $\{a_{11} = 0, b_{13} = 0, c_{12} = 0, d_{12} = 0, l_{12} = 0, m_{12} = 0\}$, $\{a_{13} = 0, b_{13} = 0, c_{12} = 0, d_{12} = 0, l_{12} = 0, m_{12} = 0\}$, $\{a_{11} = 0, a_{13} = 0, c_{12} = 0, d_{12} = 0, l_{12} = 0, m_{12} = 0\}$, $\{a_{11} = 0, a_{13} = 0, b_{13} = 0, c_{12} = 0, d_{12} = 0, l_{12} = 0\}$.  In any case, we obtain that $X\in\fr{m}_1$ or $X\in\fr{m}_2$.
}
\end{example} 


\medskip
\noindent 
{\bf Acknowledgement.}  The author would like to express her gratitude to the referee for several constructive comments on the paper.

\end{document}